\documentclass[12pt]{article} 
\usepackage{amsmath,amssymb}

\newcommand{\dd}{{\rm d}}
\newcommand{\ha}{\frac{1}{2}}
\newcommand{\mha}{-\frac{1}{2}}
\newcommand{\tha}{\frac{3}{2}}

\begin{document} 
%{ \pagestyle{empty} 

\vspace{15mm}  

%\centerline{\Large \bf  }
\centerline{\Large \bf  Jacobi's Inversion Problem}
\centerline{\Large \bf  for Genus Two Hyperelliptic Integral II}
\vspace{10mm}

\centerline{Kazuyasu Shigemoto\footnote{E-mail address:
shigemot@tezukayama-u.ac.jp}}

\centerline {{\it Tezukayama University, Tezukayama 7, Nara 631, Japan}}

\vspace{10mm}

\centerline{\bf Abstract}  
\noindent 
 In the previous paper, we reviewed the Rosenhain's paper
to the Jacobi's inversion problem for the genus two 
hyperelliptic integral.
In this paper, we review the G\"{o}pel's paper 
to the Jacobi's inversion problem for the genus two 
hyperelliptic integral.

%\newpage
%}

\vspace{20mm}
%%%%%%%%%%%%%%%%%%%%%%%%%%%%%%%%%%%%%%%%%%%%%%%%%%%%
%%%%%%%%%%%%%%%%%%%%%%% Section 1 %%%%%%%%%%%%%%%%%%%%%%%%
\section{Introduction}
\setcounter{equation}{0}
The doubly periodic function, the elliptic function, is obtained as the 
ratio of the one variable theta function. While, according to Abel, the elliptic 
function is obtained as the solution of the 
inversion problem of the elliptic integral.
Similarly, the multiply periodic function such as the genus two hyperelliptic 
function,
is obtained as the ratio of the two variables theta function.
If we obtain such multiply periodic function as the inversion problem 
of the genus two hyperelliptic integral, Jacobi noticed that 
the naive solution of the inversion problem provides the non-single valued 
multiply periodic function~\cite{Jacobi1}.
Jacobi found that the symmetric 
combination of the hyperelliptic integral becomes the single 
valued function.
Thus, Jacobi proposed the following Jacobi's inversion problem
\begin{equation}
u =\int^x \frac{\dd x}{\sqrt{f_5(x)}}+\int^{x'} \frac{\dd x'}{\sqrt{f_5(x')}},
\quad 
v =\int^x \frac{x \dd x}{\sqrt{f_5(x)}}+\int ^{x'} \frac{x' \dd x'}{\sqrt{f_5(x')}},
\nonumber
\end{equation}
for the fifth-degree polynomial function $f_5(x)$.
For this inversion problem, Jacobi conjectured that the multiply periodic 
function, which is the solution of the Jacobi's inversion problem,  is given by 
the symmetric combination such as $x+x'$, $x x'$, which become the single
valued functions of $u$ and $v$; furthermore, $x+x'$, $x x'$ are expressed 
as the ratio of the two variable theta functions~\cite{Jacobi2}.

The Jacobi's inversion problem for the genus two case is solved 
by G\"{o}pel~\cite{Gopel1,Gopel2} and independently by 
Rosenhain~\cite{Rosenhain1,Rosenhain2}.

In Rosenhain's approach, the Riemann type addition formula of 
the hyperelliptic theta function is used.
The key identities are the three quadratic theta identities, which parametrize 
five ratios of the square of hyperelliptic 
functions, which naturally gives the fifth-degree polynomial function $f_5(x)$.
While, in G\"{o}pel's approach, the theta formula by the 
duplication method of the hyperelliptic theta function is used. 
The key identities are three quartic theta identities. 
One of these identities is the Kummer's quartic 
identity, i.e., Kummer surface relation.

In the previous paper, we reviewed the Rosenhain' s paper~\cite{Shigemoto1}. In this paper, we 
review G\"{o}pel's paper.
 
%%%%%%%%%%%%%%%%%%%%%%%%%%%%%%%%%%%%%%%%%%%%%%%%%%%%%%%%%%%%%%%%%%
%%%%%%%%%%%%%%%%%%%%%%%% Section 2 %%%%%%%%%%%%%%%%%%%%%%%%%%%%%
%%%%%%%%%%%%%%%%%%%%%%%%%%%%%%%%%%%%%%%%%%%%%%%%%%%%%%%%%%%%
\section{The addition formulae of the genus two hyperelliptic theta functions
-Duplication method-}
\setcounter{equation}{0}

The quadratic addition formula, which is called the 
duplication method, was first introduced by Jacobi~\cite{Jacobi3} in 
1828 for the 
the genus one elliptic theta functions.
G\"{o}pel~\cite{Gopel1,Gopel2} used this addition formula
for the genus two hyperelliptic theta functions in 1847. 
Using this addition formula, we provided the full 
addition formula 
for the genus two hyperelliptic functions instead of the hyperelliptic
theta functions~\cite{Shigemoto2}.

The quartic addition formula, which is called the Riemann's theta formula, 
was first introduced by Jacobi~\cite{Jacobi4} in 1838 for the genus 
one elliptic theta functions.
Rosenhain~\cite{Rosenhain1, Rosenhain2} used this addition formula 
for the genus two hyperelliptic theta 
functions in 1850.  Prym~\cite{Prym} named such addition formula 
as { \it the Riemann's theta formula} because he writes down the proof 
according to the 
direct suggestion by Riemann~\cite{Riemann}  at Pisa in 1865.
Using this addition formula, the full addition formula 
for the genus two hyperelliptic functions instead of the hyperelliptic
theta functions are given by Kossak~\cite{Kossak,Shigemoto3}.

%%%%%%%%%%%%%%%%%%%%%%%%%%%%%%%%%%%%%%%%%%%%%%%%%

The theta functions with two variables are defined by
\begin{eqnarray}
&&\hskip -15mm  \vartheta\left[\begin{array}{cc} a \ c \\ b \ d \\ 
\end{array}\right](x,y;\tau_{11}, \tau_{22}, \tau_{12})
\nonumber\\
&&=\sum_{m,n\in Z} \exp 
\left\{ \pi i  
\left( \tau_{11} (m+\frac{a}{2})^2+\tau_{22} (n+\frac{c}{2})^2+
2 \tau_{12} (m+\frac{a}{2}) (n+\frac{c}{2}) \right) \right.
\nonumber\\
&& \left. +2\pi i \ 
\left( (m+\frac{a}{2})(x+\frac{b}{2})+(n+\frac{c}{2})(y+\frac{d}{2})
\right)
\right\}  ,
\label{2e1}
\end{eqnarray}
where we assume that ${\rm Im} \tau_{11}>0$, ${\rm Im}\tau_{22}>0$,
$({\rm Im} \tau_{11})({\rm Im}\tau_{22})-({\rm Im} \tau_{12})^2>0$ in order
that the summation of $m,n \in Z$ becomes convergent. 
Renaming $m\rightarrow m$, $n\rightarrow -n$, we can always choose
${\rm Im} \tau_{12}>0$, thus we assume ${\rm Im} \tau_{12}>0$.\\
We simply denote
 $\displaystyle{ \vartheta\left[\begin{array}{cc} a \ c \\ b \ d \\ \end{array}\right](x,y)
=\vartheta\left[\begin{array}{cc} a \ c \\ b \ d \\ 
\end{array}\right](x,y;\tau_{11}, \tau_{22}, \tau_{12})}$ .\\
We also use the notation  
 $\displaystyle{ \varTheta\left[\begin{array}{cc} a \ c \\ b \ d \\ \end{array}\right](x,y)
=\vartheta\left[\begin{array}{cc} a \ c \\ b \ d \\ 
\end{array}\right](x,y;2\tau_{11}, 2\tau_{22}, 2\tau_{12})}$ . 
 %According to G\"{o}pel, we use the simplified notation
%
These have the following properties 
\begin{eqnarray}
&&\hskip -5mm \vartheta\left[\begin{array}{cc} a+2 \ c \\ b\ \ \ \ d \\ \end{array}\right](x,y)
=\vartheta\left[\begin{array}{cc} a \ c+2 \\ b\ \ \ \ d \\ \end{array}\right](x,y)
=\vartheta\left[\begin{array}{cc} a \ c \\ b \ d \\ \end{array}\right](x,y),
\nonumber\\
&&\hskip -5mm \vartheta\left[\begin{array}{cc} a\ \ \ \ c \\ b+2 \ d \\ \end{array}\right](x,y)
=(\sqrt{-1})^{2a} \vartheta\left[\begin{array}{cc} a \ c \\ b \ d \\ \end{array}\right](x,y), \ 
\vartheta\left[\begin{array}{cc} a\ \ \ \ c \\ b \ d+2 \\ \end{array}\right](x,y) 
=(\sqrt{-1})^{2c} \vartheta\left[\begin{array}{cc} a \ c \\ b \ d \\ \end{array}\right](x,y) .
\nonumber
\end{eqnarray}
In order to obtain the addition formula according to the duplication method, 
we consider the product of the $\vartheta$ function in the form
\begin{eqnarray}
&&\vartheta\left[\begin{array}{cc} a \ c \\ b \ d \\ \end{array}\right](x_1,y_1) \ 
\vartheta\left[\begin{array}{cc} e \ g \\ f \ h \\ \end{array}\right](x_2,y_2) 
\nonumber\\
&&=\sum_{m,n,p,q \in Z} \exp 
\left\{ \pi i  
\left( \tau_{11} (m+\frac{a}{2})^2+\tau_{22} (n+\frac{c}{2})^2+
2 \tau_{12} (m+\frac{a}{2}) (n+\frac{c}{2}) 
\right.  \right.
\nonumber\\
&&\left. \left. +\tau_{11} (p+\frac{e}{2})^2+\tau_{22} (q+\frac{g}{2})^2+
2 \tau_{12} (p+\frac{e}{2}) (q+\frac{g}{2}) 
\right) 
\right.
\nonumber\\
&& \left. +2\pi i \ 
\left( (m+\frac{a}{2})(x_1+\frac{b}{2})+(n+\frac{c}{2})(y_1+\frac{d}{2})
  +(p+\frac{e}{2})(x_2+\frac{f}{2})+(q+\frac{g}{2})(y_2+\frac{h}{2})
\right)
\right\}  ,
\nonumber\\
&&=\sum_{m_1,n_1,m_2,n_2 \in Z} \exp 
\left\{ \frac{\pi i}{2}  
\left( \tau_{11} (m_1+\frac{a+e}{2})^2+\tau_{22} (n_1+\frac{c+g}{2})^2+
2 \tau_{12} (m_1+\frac{a+e}{2}) (n_1+\frac{c+g}{2}) 
\right.  \right.
\nonumber\\
&&\left. \left. +\tau_{11} (m_2+\frac{a-e}{2})^2+\tau_{22} (n_2+\frac{c-g}{2})^2+
2 \tau_{12} (m_2+\frac{a-e}{2}) (n_2+\frac{c-g}{2}) 
\right) 
\right.
\nonumber\\
&& \left. +\pi i \ 
\left( (m_1+\frac{a+e}{2})(x_1+x_2+\frac{b+f}{2})+(n_1+\frac{c+g}{2})(y_1+y_2+\frac{d+h}{2})
\right. \right. \nonumber\\
&& \left. \left. +(m_2+\frac{a-e}{2})(x_1-x_2+\frac{b-f}{2})+(n_2+\frac{c-g}{2})(y_1-y_2+\frac{d-h}{2})
\right)
\right\}  ,
\nonumber
\end{eqnarray}
where we use $m_1=m+p, m_2=m-p, n_1=n+q, n_2=n-q$. 
With $m_1-m_2=2p=\text{(even number)}$, $\{m_1, m_2\}$ are 
both even number
or odd number. The pair $\{n_1, n_2\}$ are also both even number or odd number.
Hence, we obtain 4 cases
i) $m_1,m_2:\text{even number}; n_1,n_2:\text{even number}$,
ii) $m_1,m_2:\text{even number}; n_1,n_2:\text{odd number}$,
iii) $m_1,m_2:\text{odd number}; n_1,n_2:\text{even number}$,
iv) $m_1,m_2:\text{odd number}; n_1,n_2:\text{odd number}$.
Therefore, we obtain 
\begin{eqnarray}
&&\vartheta\left[\begin{array}{cc} a \ c \\ b \ d \\ \end{array}\right](x_1,y_1) \ 
\vartheta\left[\begin{array}{cc} e \ g \\ f \ h \\ \end{array}\right](x_2,y_2) 
\nonumber\\
&&
=\varTheta\left[\begin{array}{cc} \frac{a+e}{2} \ \frac{c+g}{2} \\ b+f \ d+h \\ \end{array}\right](x_1+x_2 , y_1+y_2) \ \
\varTheta\left[\begin{array}{cc} \frac{a-e}{2} \ \frac{c-g}{2} \\ b-f \ d-h \\ \end{array}\right](x_1-x_2 , y_1-y_2) 
\nonumber\\
&&+\varTheta\left[\begin{array}{cc} \frac{a+e}{2} \ \frac{c+g+2}{2} \\ b+f \ d+h \\ \end{array}\right](x_1+x_2 , y_1+y_2) \ \ 
\varTheta\left[\begin{array}{cc} \frac{a-e}{2} \ \frac{c-g+2}{2} \\ b-f \ d-h \\ \end{array}\right](x_1-x_2 , y_1-y_2) \
\nonumber\\
&&+\varTheta\left[\begin{array}{cc} \frac{a+e+2}{2} \ \frac{c+g}{2} \\ b+f \ d+h \\ \end{array}\right](x_1+x_2 , y_1+y_2) \ \
\varTheta\left[\begin{array}{cc} \frac{a-e+2}{2} \ \frac{c-g}{2} \\ b-f \ d-h \\ \end{array}\right](x_1-x_2 , y_1-y_2) 
\nonumber\\
&&+\varTheta\left[\begin{array}{cc} \frac{a+e+2}{2} \ \frac{c+g+2}{2} \\ b+f \ d+h \\ \end{array}\right](x_1+x_2 , y_1+y_2) \ \
\varTheta\left[\begin{array}{cc} \frac{a-e+2}{2} \ \frac{c-g+2}{2} \\ b-f \ d-h \\ \end{array}\right](x_1-x_2 , y_1-y_2) .
\nonumber
\end{eqnarray}
Putting $x_1=u_1+u_2$, $x_2=u_1-u_2$, $y_1=v_1+v_2$, $y_2=v_1-v_2$, we express
the above formula
in the form
\begin{eqnarray}
&&\vartheta\left[\begin{array}{cc} a \ c \\ b \ d \\ \end{array}\right](u_1+u_2 ,v_1+v_2) \ 
\vartheta\left[\begin{array}{cc} e \ g \\ f \ h \\ \end{array}\right](u_1-u_2, v_1-v_2) 
\nonumber\\
&&
=\varTheta\left[\begin{array}{cc} \frac{a+e}{2} \ \frac{c+g}{2} \\ b+f \ d+h \\ \end{array}\right](2u_1 , 2v_1) \ \ 
\varTheta\left[\begin{array}{cc} \frac{a-e}{2} \ \frac{c-g}{2} \\ b-f \ d-h \\ \end{array}\right](2u_2 , 2v_2) 
\nonumber\\
&&+\varTheta\left[\begin{array}{cc} \frac{a+e}{2} \ \frac{c+g+2}{2} \\ b+f \ d+h \\ \end{array}\right](2u_1 , 2v_1) \ \
\varTheta\left[\begin{array}{cc} \frac{a-e}{2} \ \frac{c-g+2}{2} \\ b-f \ d-h \\ \end{array}\right](2u_2 , 2v_2) \
\nonumber\\
&&+\varTheta\left[\begin{array}{cc} \frac{a+e+2}{2} \ \frac{c+g}{2} \\ b+f \ d+h \\ \end{array}\right](2u_1 , 2v_1) \ \
\varTheta\left[\begin{array}{cc} \frac{a-e+2}{2} \ \frac{c-g}{2} \\ b-f \ d-h \\ \end{array}\right](2u_2 , 2v_2) 
\nonumber\\
&&+\varTheta\left[\begin{array}{cc} \frac{a+e+2}{2} \ \frac{c+g+2}{2} \\ b+f \ d+h \\ \end{array}\right](2u_1 , 2v_1) \ \
\varTheta\left[\begin{array}{cc} \frac{a-e+2}{2} \ \frac{c-g+2}{2} \\ b-f \ d-h \\ \end{array}\right](2u_2 , 2v_2) .
\label{2e2}
\end{eqnarray}
In the special case of $u_2=0, v_2=0$, we obtain
\begin{eqnarray}
&&\vartheta\left[\begin{array}{cc} a \ c \\ b \ d \\ \end{array}\right](u_1 ,v_1) \ 
\vartheta\left[\begin{array}{cc} e \ g \\ f \ h \\ \end{array}\right](u_1, v_1) 
\nonumber\\
&&
=\varTheta\left[\begin{array}{cc} \frac{a+e}{2} \ \frac{c+g}{2} \\ b+f \ d+h \\ \end{array}\right](2u_1 , 2v_1) \ \ 
\varTheta\left[\begin{array}{cc} \frac{a-e}{2} \ \frac{c-g}{2} \\ b-f \ d-h \\ \end{array}\right](0, 0) 
\nonumber\\
&&+\varTheta\left[\begin{array}{cc} \frac{a+e}{2} \ \frac{c+g+2}{2} \\ b+f \ d+h \\ \end{array}\right](2u_1 , 2v_1) \ \
\varTheta\left[\begin{array}{cc} \frac{a-e}{2} \ \frac{c-g+2}{2} \\ b-f \ d-h \\ \end{array}\right](0 , 0) \
\nonumber\\
&&+\varTheta\left[\begin{array}{cc} \frac{a+e+2}{2} \ \frac{c+g}{2} \\ b+f \ d+h \\ \end{array}\right](2u_1 , 2v_1) \ \
\varTheta\left[\begin{array}{cc} \frac{a-e+2}{2} \ \frac{c-g}{2} \\ b-f \ d-h \\ \end{array}\right](0 , 0) 
\nonumber\\
&&+\varTheta\left[\begin{array}{cc} \frac{a+e+2}{2} \ \frac{c+g+2}{2} \\ b+f \ d+h \\ \end{array}\right](2u_1 , 2v_1) \ \
\varTheta\left[\begin{array}{cc} \frac{a-e+2}{2} \ \frac{c-g+2}{2} \\ b-f \ d-h \\ \end{array}\right](0 , 0) .
\label{2e3}
\end{eqnarray}
We frequently use in the following form 
\begin{eqnarray}
  &&\hspace{-20mm} 
  \left( \begin{array}{c}
  \vartheta[\begin{array}{cc} a \ b \\ 0 \ 0 \\ \end{array}]^2(u ,v) \\
  \vartheta[\begin{array}{cc} a \ b \\ 0 \ 1 \\ \end{array}]^2(u ,v) \\
  \vartheta[\begin{array}{cc} a \ b \\ 1 \ 0 \\ \end{array}]^2(u ,v) \\
  \vartheta[\begin{array}{cc} a \ b \\ 1 \ 1 \\ \end{array}]^2(u ,v) 
 \end{array} \right)
 =\left( \begin{array}{cccc} 
  1 & 1& 1& 1\\ 
  1 & -1& 1 & -1\\
  1 & 1& -1 & -1\\ 
  1 & -1& -1 & 1\end{array} \right)
  \left( \begin{array}{c}
  \varTheta[\begin{array}{cc} 0 \ 0 \\ 0 \ 0 \\ \end{array}](2u ,2v) 
  \varTheta[\begin{array}{cc} a \ b \\ 0 \ 0 \\ \end{array}](0 ,0) \\
  \varTheta[\begin{array}{cc} 0 \ 1 \\ 0 \ 0 \\ \end{array}](2u ,2v) 
  \varTheta[\begin{array}{cc} a \ b-1 \\ 0 \ 0 \\ \end{array}](0 ,0)\\
  \varTheta[\begin{array}{cc} 1 \ 0 \\ 0 \ 0 \\ \end{array}](2u ,2v)  
  \varTheta[\begin{array}{cc} a-1 \ b \\ 0 \ 0 \\ \end{array}](0 ,0)  \\
  \varTheta[\begin{array}{cc} 1 \ 1 \\ 0 \ 0 \\ \end{array}](2u ,2v) 
  \varTheta[\begin{array}{cc} a-1 \ b-1 \\ 0 \ 0 \\ \end{array}](0 ,0) 
 \end{array} \right)   .
\nonumber
\end{eqnarray}
Here the Riemann matrix, which is the self-adjoint orthogonal matrix, emerges. 
We inversely solve the above in the following matrix form
\begin{eqnarray}
  &&\hspace{-20mm} 
  \left( \begin{array}{c}
  \varTheta[\begin{array}{cc} 0 \ 0 \\ 0 \ 0 \\ \end{array}](2u ,2v) 
  \varTheta[\begin{array}{cc} a \ b \\ 0 \ 0 \\ \end{array}](0 ,0) \\
  \varTheta[\begin{array}{cc} 0 \ 1 \\ 0 \ 0 \\ \end{array}](2u ,2v) 
  \varTheta[\begin{array}{cc} a \ b-1 \\ 0 \ 0 \\ \end{array}](0 ,0)\\
  \varTheta[\begin{array}{cc} 1 \ 0 \\ 0 \ 0 \\ \end{array}](2u ,2v)  
  \varTheta[\begin{array}{cc} a-1 \ b \\ 0 \ 0 \\ \end{array}](0 ,0)  \\
  \varTheta[\begin{array}{cc} 1 \ 1 \\ 0 \ 0 \\ \end{array}](2u ,2v) 
  \varTheta[\begin{array}{cc} a-1 \ b-1 \\ 0 \ 0 \\ \end{array}](0 ,0) 
  \end{array} \right)
 =\frac{1}{4}\left( \begin{array}{cccc} 
  1 & 1& 1& 1\\ 
  1 & -1& 1 & -1\\
  1 & 1& -1 & -1\\ 
  1 & -1& -1 & 1\end{array} \right)
  \left( \begin{array}{c}
  \vartheta[\begin{array}{cc} a \ b \\ 0 \ 0 \\ \end{array}]^2(u ,v) \\
  \vartheta[\begin{array}{cc} a \ b \\ 0 \ 1 \\ \end{array}]^2(u ,v) \\
  \vartheta[\begin{array}{cc} a \ b \\ 1 \ 0 \\ \end{array}]^2(u ,v) \\
  \vartheta[\begin{array}{cc} a \ b \\ 1 \ 1 \\ \end{array}]^2(u ,v) 
 \end{array} \right)   . 
\nonumber
\end{eqnarray}
%
%%%%%%%%%%%%%%%%%%%%%%%%%%%%%%%%%%
We use the G\"{o}pel's simplified notation of the form
\begin{eqnarray}
&&\hskip -20mm P_0=\vartheta[\begin{array}{cc} 0 \ 0 \\ 1 \ 1 \\ \end{array}](u ,v), \quad
P_1=\vartheta[\begin{array}{cc} 0 \ 0 \\ 0 \ 1 \\ \end{array}](u ,v), \quad
P_2=\vartheta[\begin{array}{cc} 0 \ 0 \\ 1 \ 0 \\ \end{array}](u ,v), \quad
P_3=\vartheta[\begin{array}{cc} 0 \ 0 \\ 0 \ 0 \\ \end{array}](u ,v), 
\nonumber\\
&&\hskip -20mm Q_0=\vartheta[\begin{array}{cc} 1 \ 0 \\ 1 \ 1 \\ \end{array}](u ,v), \quad
Q_1=\vartheta[\begin{array}{cc} 1 \ 0 \\ 0 \ 1 \\ \end{array}](u ,v), \quad
Q_2=\vartheta[\begin{array}{cc} 1 \ 0 \\ 1 \ 0 \\ \end{array}](u ,v), \quad
Q_3=\vartheta[\begin{array}{cc} 1 \ 0 \\ 0 \ 0 \\ \end{array}](u ,v), 
\nonumber\\
&&\hskip -20mm R_0=\vartheta[\begin{array}{cc} 0 \ 1 \\ 1 \ 1 \\ \end{array}](u ,v), \quad
R_1=\vartheta[\begin{array}{cc} 0 \ 1 \\ 0 \ 1 \\ \end{array}](u ,v), \quad
R_2=\vartheta[\begin{array}{cc} 0 \ 1 \\ 1 \ 0 \\ \end{array}](u ,v), \quad
R_3=\vartheta[\begin{array}{cc} 0 \ 1 \\ 0 \ 0 \\ \end{array}](u ,v), 
\nonumber\\
&&\hskip -20mm S_0=\vartheta[\begin{array}{cc} 1 \ 1 \\ 1 \ 1 \\ \end{array}](u ,v), \quad
S_1=\vartheta[\begin{array}{cc} 1 \ 1 \\ 0 \ 1 \\ \end{array}](u ,v), \quad
S_2=\vartheta[\begin{array}{cc} 1 \ 1 \\ 1 \ 0 \\ \end{array}](u ,v), \quad
S_3=\vartheta[\begin{array}{cc} 1 \ 1 \\ 0 \ 0 \\ \end{array}](u ,v), 
\nonumber\\
&&\hskip -20mm
T=\varTheta\left[\begin{array}{cc} 0 \ 0 \\ 0 \ 0 \\ \end{array}\right](2u, 2v), \
U=\varTheta\left[\begin{array}{cc} 1 \ 0 \\ 0 \ 0 \\ \end{array}\right](2u, 2v), \ 
V=\varTheta\left[\begin{array}{cc} 0 \ 1 \\ 0 \ 0 \\ \end{array}\right](2u, 2v), \ 
W=\varTheta\left[\begin{array}{cc} 1 \ 1 \\ 0 \ 0 \\ \end{array}\right](2u, 2v), 
\nonumber\\
&&\hskip -20mm
t=\varTheta\left[\begin{array}{cc} 0 \ 0 \\ 0 \ 0 \\ \end{array}\right](0, 0), \
u=\varTheta\left[\begin{array}{cc} 1 \ 0 \\ 0 \ 0 \\ \end{array}\right](0, 0), \ 
v=\varTheta\left[\begin{array}{cc} 0 \ 1 \\ 0 \ 0 \\ \end{array}\right](0, 0), \ 
w=\varTheta\left[\begin{array}{cc} 1 \ 1 \\ 0 \ 0 \\ \end{array}\right](0, 0)   .  
 \label{2e4}
\end{eqnarray}
Using the G\"{o}pel's notation, we obtain 
\begin{eqnarray}
  &&\hspace{-30mm} 
  \left( \begin{array}{c}
  P_3^2 \\ P_1^2 \\P_2^2 \\P_0^2 \end{array} \right)
 =\left( \begin{array}{cccc} 
  1 & 1& 1& 1\\ 
  1 & -1& 1 & -1\\
  1 & 1& -1 & -1\\ 
  1 & -1& -1 & 1\end{array} \right)
  \left( \begin{array}{c}
  t T \\ v V \\ u U \\ w W \end{array} \right)   , \ 
  \left( \begin{array}{c}
  t T \\ v V \\ u U \\ w W  \\P_0^2 \end{array} \right)
 =\frac{1}{4}\left( \begin{array}{cccc} 
  1 & 1& 1& 1\\ 
  1 & -1& 1 & -1\\
  1 & 1& -1 & -1\\ 
  1 & -1& -1 & 1\end{array} \right)
  \left( \begin{array}{c}P_3^2 \\ P_1^2 \\P_2^2 \\P_0^2 \end{array} \right) , 
 \label{2e5}
\end{eqnarray}
\begin{eqnarray}
  &&\hspace{-30mm} 
  \left( \begin{array}{c}
  Q_3^2 \\ Q_1^2 \\Q_2^2 \\Q_0^2 \end{array} \right)
 =\left( \begin{array}{cccc} 
  1 & 1& 1& 1\\ 
  1 & -1& 1 & -1\\
  1 & 1& -1 & -1\\ 
  1 & -1& -1 & 1\end{array} \right)
  \left( \begin{array}{c}
  u T \\ w V \\ t U \\ v W \end{array} \right)   ,\ 
  \left( \begin{array}{c}
  u T \\ w V \\ t U \\ v W 2 \end{array} \right)
 =\frac{1}{4}\left( \begin{array}{cccc} 
  1 & 1& 1& 1\\ 
  1 & -1& 1 & -1\\
  1 & 1& -1 & -1\\ 
  1 & -1& -1 & 1\end{array} \right)
  \left( \begin{array}{c}
 Q_3^2 \\ Q_1^2 \\Q_2^2 \\Q_0^2 \end{array} \right) ,
 \label{2e6}
\end{eqnarray}
\begin{eqnarray}
  &&\hspace{-30mm} 
  \left( \begin{array}{c}
  R_3^2 \\ R_1^2 \\R_2^2 \\R_0^2 \end{array} \right)
 =\left( \begin{array}{cccc} 
  1 & 1& 1& 1\\ 
  1 & -1& 1 & -1\\
  1 & 1& -1 & -1\\ 
  1 & -1& -1 & 1\end{array} \right)
  \left( \begin{array}{c}
  v T \\ t V \\ w U \\ u W \end{array} \right)   ,\ 
  \left( \begin{array}{c}
  v T \\ t V \\ w U \\ u W \end{array} \right)
 =\frac{1}{4}\left( \begin{array}{cccc} 
  1 & 1& 1& 1\\ 
  1 & -1& 1 & -1\\
  1 & 1& -1 & -1\\ 
  1 & -1& -1 & 1\end{array} \right)
  \left( \begin{array}{c}
  R_3^2 \\ R_1^2 \\R_2^2 \\R_0^2 \end{array} \right) ,
 \label{2e7}
\end{eqnarray}
\begin{eqnarray}
  &&\hspace{-30mm} 
  \left( \begin{array}{c}
  S_3^2 \\ S_1^2 \\S_2^2 \\S_0^2 \end{array} \right)
 =\left( \begin{array}{cccc} 
  1 & 1& 1& 1\\ 
  1 & -1& 1 & -1\\
  1 & 1& -1 & -1\\ 
  1 & -1& -1 & 1\end{array} \right)
  \left( \begin{array}{c}
  w T \\ u V \\ v U \\ t W \end{array} \right)   ,\ 
\left( \begin{array}{c}
  w T \\ u V \\ v U \\ t W \end{array} \right)
 =\frac{1}{4}\left( \begin{array}{cccc} 
  1 & 1& 1& 1\\ 
  1 & -1& 1 & -1\\
  1 & 1& -1 & -1\\ 
  1 & -1& -1 & 1\end{array} \right)
  \left( \begin{array}{c}
  S_3^2 \\ S_1^2 \\S_2^2 \\S_0^2 \end{array} \right) .
 \label{2e8}
\end{eqnarray}
%

%%%%%%%%%%%%%%%%%%%%%%%%%%%%%%%%%%%%%%%%%%%%%%%%%%%%%%%%%%%%%%%%%%%%%
%%%%%%%%%%%%%%%%%%%%%%%%% Section 3 %%%%%%%%%%%%%%%%%%%%%%%%%%%%%%%%%%% 
%%%%%%%%%%%%%%%%%%%%%%%%%%%%%%%%%%%%%%%%%%%%%%%%%%%%%%%%%%%%
\section{The Kummer surface relation}
\setcounter{equation}{0}
We use the Kummer's quartic relation, i.e., Kummer surface relation, 
 to solve the Jacobi's inversion relation. Thus,  we
derive the Kummer's quartic relation.
By using the addition formula of $P_1$, $S_1$, $P_2$ and $S_2$, we obtain 
\begin{eqnarray}
&&P_1^2=tT+uU-vV-wW, \quad P_2^2=tT-uU+vV-wW, 
\label{3e1}\\
&&S_1^2=wT+vU-uV-tW, \quad S_2^2=wT-vU+uV-tW . 
\label{3e2}
\end{eqnarray}
We inversely solve in the form
\begin{eqnarray}
&&\hskip -20mm T=\frac{ t (P_1^2+P_2^2) -w (S_1^2+S_2^2)}{2(t^2-w^2)} , \quad
U=\frac{ u (P_1^2-P_2^2) -v (S_1^2-S_2^2)}{2(u^2-v^2)} ,
\label{3e3}\\
&&\hskip -20mm V=\frac{v (P_1^2-P_2^2) -u (S_1^2-S_2^2)}{2(u^2-v^2)} ,  \quad
W=\frac{w (P_1^2+P_2^2) -t (S_1^2+S_2^2)}{2(t^2-w^2)} .
\label{3e4}
\end{eqnarray}
Later, we use  the following
\begin{eqnarray}
&&Q_1^2=uT+tU-wV-vW, \quad R_1^2=vT+wU-tV-uW . 
\label{3e5}
\end{eqnarray}
%

%%%%%%%%%%%%%%%%%%%%%%%%%%%%%%%%%%%%%
\subsection{Proof of $Q_1 R_1=a P_1 S_1+bP_2 S_2$}
We will demonstrate that $Q_1 R_1=a P_1 S_1 +b P_2 S_2$, $(a,b =\text{const.})$. Considering the square of this relation, we obtain
$Q_1^2 R_1^2=(a P_1 S_1 +b P_2 S_2)^2$. Using  Eqs.(\ref{3e3})-(\ref{3e5}), the left-hand 
side term $Q_1^2 R_1^2$ is expressed 
as the polynomial of $P_1, S_1, P_2, S_2$.
While the right-hand side term $(a P_1 S_1 +b P_2 S_2)^2$ is of course the polynomial of $P_1, S_1, P_2, S_2$. Hence, $Q_1^2 R_1^2=(a P_1 S_1 +b P_2 S_2)^2$ gives the Kummer's quartic relations of $P_1, S_1, P_2, S_2$.
If we can find one of the expressions, which is expressed as the 
linear combination of $P_1 S_1$ and $P_2 S_2$, it suffices for our purpose. For example, 
instead of using $Q_1 R_1$, we can use $Q_2 R_2$ because this also becomes the 
linear combination of  $P_1 S_1$ and $P_2 S_2$.
By using the addition formula, we obtain 
\begin{eqnarray}
&&\hskip -10mm Q_1 R_1=\vartheta[\begin{array}{cc} 1 \ 0 \\ 0 \ 1 \\ \end{array}](u ,v)
\vartheta[\begin{array}{cc} 0 \ 1 \\ 0 \ 1 \\ \end{array}](u ,v)
\nonumber\\
&&\hskip -10mm =\varTheta[\begin{array}{cc} \ha \ \ha \\ 0\ \ 2 \\ \end{array}](2u, 2v) 
\varTheta[\begin{array}{cc} \ha \ \mha \\ 0\ \ 0 \\ \end{array}](0, 0)
%\nonumber\\
%&&
+\varTheta[\begin{array}{cc} \ha \ \tha \\ 0\ \ 2 \\ \end{array}](2u, 2v) 
\varTheta[\begin{array}{cc} \ha \ \ha \\ 0\ \ 0 \\ \end{array}](0, 0)
\nonumber\\
&&\hskip -10mm +\varTheta[\begin{array}{cc} \tha \ \ha \\ 0\ \ 2 \\ \end{array}](2u, 2v) 
\varTheta[\begin{array}{cc} \tha \ \mha \\ 0\ \ \ 0 \\ \end{array}](0, 0)
%\nonumber\\
%&&
+\varTheta[\begin{array}{cc}  \tha \ \tha \\ 0\ \ 2 \\ \end{array}](2u, 2v) 
\varTheta[\begin{array}{cc} \tha \ \ha \\ 0\ \ 0 \\ \end{array}](0, 0)  ,
%\label{3e6}\\
\nonumber
\end{eqnarray}
which provides
\begin{eqnarray}
&&\hskip -10mm Q_1 R_1=i \varTheta[\begin{array}{cc} \ha \ \ha \\ 0\ \ 0 \\ \end{array}](2u, 2v) 
\varTheta[\begin{array}{cc} \ha \ \mha \\ 0\ \ \ 0 \\ \end{array}](0, 0)
%\nonumber\\
%&&
-i\varTheta[\begin{array}{cc} \ha \ \mha \\ 0\ \ 0 \\ \end{array}](2u, 2v) 
\varTheta[\begin{array}{cc} \ha \ \ha \\ 0\ \ 0 \\ \end{array}](0, 0)
\nonumber\\
&&\hskip -10mm 
+i \varTheta[\begin{array}{cc} \mha \ \ha \\ 0\ \ 0 \\ \end{array}](2u, 2v) 
\varTheta[\begin{array}{cc} \mha \ \mha \\ 0\ \ \ 0 \\ \end{array}](0, 0)
%\nonumber\\
%&&
-i\varTheta[\begin{array}{cc} \mha \ \mha \\ 0\ \ 0 \\ \end{array}](2u, 2v) 
\varTheta[\begin{array}{cc} \mha \ \ha \\ 0\ \ 0 \\ \end{array}](0, 0)
\nonumber\\
&&\hskip -10mm 
=i\left(\varTheta[\begin{array}{cc} \ha \ \ha \\ 0\ \ 0 \\ \end{array}](2u, 2v) 
-\varTheta[\begin{array}{cc} \mha \ \mha \\ 0\ \ 0 \\ \end{array}](2u, 2v)\right)
\varTheta[\begin{array}{cc} \ha \mha \\ 0\ \ \ 0 \\ \end{array}](0, 0)
\nonumber\\
&&\hskip -10mm 
-i\left(\varTheta[\begin{array}{cc} \ha \ \mha \\ 0\ \ 0 \\ \end{array}](2u, 2v) 
-\varTheta[\begin{array}{cc} \mha \ \ha \\ 0\ \ 0 \\ \end{array}](2u, 2v)\right)
\varTheta[\begin{array}{cc} \ha \ \ha \\ 0\ \ \ 0 \\ \end{array}](0, 0)  ,
\label{3e6}
\end{eqnarray}
where we used
\begin{eqnarray}
&&\varTheta[\begin{array}{cc} \ha \ \ha \\ 0\ \ \ 0 \\ \end{array}](0, 0)
=\varTheta[\begin{array}{cc} \mha \ \mha \\ 0\ \ \ 0 \\ \end{array}](0, 0), \quad
%\nonumber\\
%&&
\varTheta[\begin{array}{cc} \ha \ \mha \\ 0\ \ \ 0 \\ \end{array}](0, 0)
=\varTheta[\begin{array}{cc} \mha \ \ha \\ 0\ \ \ 0 \\ \end{array}](0, 0)  . 
\label{3e7}
\end{eqnarray}
Similarly,  we obtain
\begin{eqnarray}
&&\hskip -10mm P_1 S_1=\vartheta[\begin{array}{cc} 0 \ 0 \\ 0 \ 1 \\ \end{array}](u ,v)
\vartheta[\begin{array}{cc} 1 \ 1 \\ 0 \ 1 \\ \end{array}](u ,v)
\nonumber\\
&&\hskip -10mm 
=i\left(\varTheta[\begin{array}{cc} \ha \ \ha \\ 0\ \ 0 \\ \end{array}](2u, 2v) 
-\varTheta[\begin{array}{cc} \mha \ \mha \\ 0\ \ 0 \\ \end{array}](2u, 2v)\right)
\varTheta[\begin{array}{cc} \ha \ \ha \\ 0\ \ \ 0 \\ \end{array}](0, 0) 
\nonumber\\
&&\hskip -10mm 
-i\left(\varTheta[\begin{array}{cc} \ha \ \mha \\ 0\ \ 0 \\ \end{array}](2u, 2v) 
-\varTheta[\begin{array}{cc} \mha \ \ha \\ 0\ \ 0 \\ \end{array}](2u, 2v)\right)
\varTheta[\begin{array}{cc} \ha \ \mha \\ 0\ \ \ 0 \\ \end{array}](0, 0), 
\label{3e8}\\
&&\hskip -10mm P_2 S_2=\vartheta[\begin{array}{cc} 0 \ 0 \\ 1 \ 0 \\ \end{array}](u ,v)
\vartheta[\begin{array}{cc} 1 \ 1 \\ 1 \ 0 \\ \end{array}](u ,v)
\nonumber\\
&&\hskip -10mm 
=i\left(\varTheta[\begin{array}{cc} \ha \ \ha \\ 0\ \ 0 \\ \end{array}](2u, 2v) 
-\varTheta[\begin{array}{cc} \mha \ \mha \\ 0\ \ 0 \\ \end{array}](2u, 2v)\right)
\varTheta[\begin{array}{cc} \ha \ \ha \\ 0\ \ \ 0 \\ \end{array}](0, 0)
\nonumber\\
&&\hskip -10mm 
+i\left(\varTheta[\begin{array}{cc} \ha \ \mha \\ 0\ \ 0 \\ \end{array}](2u, 2v) 
-\varTheta[\begin{array}{cc} \mha \ \ha \\ 0\ \ 0 \\ \end{array}](2u, 2v)\right)
\varTheta[\begin{array}{cc} \ha \ \mha \\ 0\ \ \ 0 \\ \end{array}](0, 0).
\label{3e9}
\end{eqnarray} 
Combining these relations, we obtain three coupled equations 
\begin{eqnarray}
&&Q_1 R_1=i B X-i A Y, 
\label{3e10}\\
&&P_1 S_1=i A X-i B Y, 
\label{3e11}\\
&&P_2 S_2=i A X+ i B Y, 
\label{3e12}\\
&&X=\left(\varTheta[\begin{array}{cc} \ha \ \ha \\ 0\ \ \ 0 \\ \end{array}](2u, 2v) 
-\varTheta[\begin{array}{cc} \mha \ \mha \\ 0\ \ \ \ 0 \ \\ \end{array}](2u, 2v)\right),
\label{3e13}\\
&&Y=\left(\varTheta[\begin{array}{cc} \ha \ \mha \\ 0 \ \ \ 0\ \ \\ \end{array}](2u, 2v) 
-\varTheta[\begin{array}{cc} \mha \ \ha \\ 0\ \ \ 0 \\ \end{array}](2u, 2v)\right) ,
\label{3e14}\\
&&A=\varTheta[\begin{array}{cc} \ha \ \ha \\ 0\ \ 0 \\ \end{array}](0, 0), \quad
B=\varTheta[\begin{array}{cc} \ha \ \mha \\ 0\ \ \ 0 \\ \end{array}](0, 0)  .
\label{3e15}
\end{eqnarray}
Expressing $X$ and $Y$ by $P_1 S_1$ and $P_2 S_2$, and substituting
 into the right-hand side of 
 $Q_1 R_1$, we complete the proof of the form 
\begin{eqnarray}
&&Q_1 R_1=a P_1 S_1+b P_2 S_2, \quad  
a=\frac{A^2+B^2}{2 A B}, \quad b=\frac{A^2-B^2}{2 A B},\quad a^2-b^2=1 .
\label{3e16}
\end{eqnarray}
%%%%%%%%%%%%%%%%%%%%%%%%%%%%%
\subsection{The expression of $a$ and $b$ with $t, u, v, w$ }
In order to express $A$ and $B$ as the function of $t, u, v, w$, we consider
the quantity 
\begin{eqnarray}
&&\hskip -10mm P_3 S_3=\vartheta[\begin{array}{cc} 1 \ 1 \\ 0\ \ 0 \\ \end{array}](u, v) 
\vartheta[\begin{array}{cc} 0 \ 0 \\ 0\ \ 0 \\ \end{array}](u, v)
\nonumber\\
&&\hskip -10mm =\left(\varTheta[\begin{array}{cc} \ha \ \ha \\ 0\ \ 0 \\ \end{array}](2u, 2v) 
+\varTheta[\begin{array}{cc} \mha \ \mha \\ 0\ \ 0 \\ \end{array}](2u, 2v)\right)
\varTheta[\begin{array}{cc} \ha \ \ha \\ 0\ \ \ 0 \\ \end{array}](0, 0)
\nonumber\\
&&\hskip -10mm 
+\left(\varTheta[\begin{array}{cc} \ha \ \mha \\ 0\ \ 0 \\ \end{array}](2u, 2v) 
+\varTheta[\begin{array}{cc} \mha \ \ha \\ 0\ \ 0 \\ \end{array}](2u, 2v)\right)
\varTheta[\begin{array}{cc} \ha \ \mha \\ 0\ \ \ 0 \\ \end{array}](0, 0) .
\label{3e17}
\end{eqnarray} 
By putting $u=0, v=0$, we obtain 
\begin{eqnarray}
&&\hskip -20mm \vartheta[\begin{array}{cc} 1 \ 1 \\ 0\ \ 0 \\ \end{array}](0, 0) 
\vartheta[\begin{array}{cc} 0 \ 0 \\ 0\ \ 0 \\ \end{array}](0, 0)
%\nonumber\\
%&&
=2\varTheta[\begin{array}{cc} \ha \ \ha \\ 0\ \ 0 \\ \end{array}](0, 0)^2 
+2\varTheta[\begin{array}{cc} \ha \ \mha \\ 0\ \ 0 \\ \end{array}](0, 0)^2
=2(A^2+B^2)  .
\label{3e18}
\end{eqnarray} 

Similarly, we consider the quantity
\begin{eqnarray}
&&\hskip -10mm Q_3 R_3=\vartheta[\begin{array}{cc} 1 \ 0 \\ 0\ \ 0 \\ \end{array}](u, v) 
\vartheta[\begin{array}{cc} 0 \ 1 \\ 0\ \ 0 \\ \end{array}](u, v)
\nonumber\\
&&\hskip -10mm =\left(\varTheta[\begin{array}{cc} \ha \ \ha \\ 0\ \ 0 \\ \end{array}](2u, 2v) 
+\varTheta[\begin{array}{cc} \mha \ \mha \\ 0\ \ 0 \\ \end{array}](2u, 2v)\right)
\varTheta[\begin{array}{cc} \ha \ \mha \\ 0\ \ \ 0 \\ \end{array}](0, 0)
\nonumber\\
&&\hskip -10mm 
+\left(\varTheta[\begin{array}{cc} \ha \ \mha \\ 0\ \ 0 \\ \end{array}](2u, 2v) 
+\varTheta[\begin{array}{cc} \mha \ \ha \\ 0\ \ 0 \\ \end{array}](2u, 2v)\right)
\varTheta[\begin{array}{cc} \ha \ \ha \\ 0\ \ \ 0 \\ \end{array}](0, 0)  .
\label{3e19}
\end{eqnarray} 
By putting $u=0, v=0$, we obtain
\begin{eqnarray}
&&\hskip -10mm \vartheta[\begin{array}{cc} 1 \ 0 \\ 0\ \ 0 \\ \end{array}](0, 0) 
\vartheta[\begin{array}{cc} 0 \ 1 \\ 0\ \ 0 \\ \end{array}](0, 0)
%\nonumber\\
%&&
=4\varTheta[\begin{array}{cc} \ha \ \ha \\ 0\ \ 0 \\ \end{array}](0, 0) 
\varTheta[\begin{array}{cc} \ha \ \mha \\ 0\ \ 0 \\ \end{array}](0, 0)
=4 AB  .
\label{3e20}
\end{eqnarray} 
Therefore, the constant $a$ is solved by the function of 
$\vartheta[\begin{array}{cc} a \ c \\ b \ d \\ \end{array}](0, 0)$ in the form
\begin{eqnarray}
a=\frac{A^2+B^2}{2AB}=\frac{\vartheta[\begin{array}{cc} 1 \ 1 \\ 0\ \ 0 \\ \end{array}](0, 0) 
\vartheta[\begin{array}{cc} 0 \ 0 \\ 0\ \ 0 \\ \end{array}](0, 0)}
{\vartheta[\begin{array}{cc} 1 \ 0 \\ 0\ \ 0 \\ \end{array}](0, 0) 
\vartheta[\begin{array}{cc} 0 \ 1 \\ 0\ \ 0 \\ \end{array}](0, 0)}   .
\label{3e21}
\end{eqnarray}
Taking the square and use the addition formula, 
$S_3(0,0)^2
=\vartheta[\begin{array}{cc} 1 \ 1 \\ 0\ \ 0 \\ \end{array}](0, 0) ^2=2(t w+u v)$,
$P_3(0,0)^2= \vartheta[\begin{array}{cc} 0 \ 0 \\ 0\ \ 0 \\ \end{array}](0, 0)^2
=t^2+u^2+v^2+w^2$, 
$Q_3(0,0)^2=\vartheta[\begin{array}{cc} 1 \ 0 \\ 0\ \ 0 \\ \end{array}](0, 0)^2 
=2(t u+v w)$, $R_3(0,0)^3= \vartheta[\begin{array}{cc} 0 \ 1 \\ 0\ \ 0 \\ \end{array}](0, 0)^2
=2(tv+u w)$, we obtain the expression of $a^2$ with $t, u, v, w$ 
\begin{eqnarray}
a^2=\frac{\vartheta[\begin{array}{cc} 1 \ 1 \\ 0\ \ 0 \\ \end{array}](0, 0) ^2
\vartheta[\begin{array}{cc} 0 \ 0 \\ 0\ \ 0 \\ \end{array}](0, 0)^2}
{\vartheta[\begin{array}{cc} 1 \ 0 \\ 0\ \ 0 \\ \end{array}](0, 0)^2 
\vartheta[\begin{array}{cc} 0 \ 1 \\ 0\ \ 0 \\ \end{array}](0, 0)^2}
=\frac{(t w + u v)(t^2+u^2+v^2+w^2)}{2(t u+v w)(t v+u w)}   .
\label{3e22}
\end{eqnarray}
By using $a^2-b^2=1$, we obtain $b^2$ in the form
\begin{eqnarray}
b^2=a^2-1=\frac{(t^2-u^2-v^2+w^2)(t w - u v)}{2(t u+v w)(t v+u w)}, 
\label{3e23}
\end{eqnarray}
which provides
\begin{eqnarray}
a=\sqrt{\frac{(t^2+u^2+v^2+w^2)(t w + u v)}{2(t u+v w)(t v+u w)}}, \quad
b=\sqrt{\frac{(t^2-u^2-v^2+w^2)(t w - u v)}{2(t u+v w)(t v+u w)}} .
\label{3e24}
\end{eqnarray}
%
%%%%%%%%%%%%%%%%%%%%%%%%%%%%%%
\subsection{The Kummer surface relation of $P_1, S_1, P_2, S_2$}
The Kummer's quartic relation
for $P_1, S_1, P_2, S_2$ is given by considering the square of 
$Q_1 R_1=a P_1 S_1+b P_2 S_2$
in the form 
$Q_1^2 R_1^2=a^2 P_1^2 S_1^2+b^2 P_2^2 S_2^2 +2 a b P_1 S_1 P_2 S_2$ 
via Eqs.(\ref{3e3})-(\ref{3e5}).
After the straightforward calculation, we obtain the following 
Kummer's quartic relation
\begin{eqnarray}
&&\hskip -10mm P_1^4+S_1^4+P_2^4+S_2^4-2F(P_1^2 P_2^2
+S_1^2 S_2^2)+2C(P_1^2 S_2^2+P_2^2 S_1^2)
\nonumber\\
&&\hskip -10mm -2E(P_1^2 S_1^2+P_2^2 S_2^2) -4D P_1 P_2 S_1 S_2=0, 
\label{3e25}\\
&&\hskip -10mm E=\frac{t^2 u^2+t^2 v^2+u^2 w^2
+v^2 w^2+ 4 t u v w}{2(t u+v w)( t v+u w)} ,
\label{3e26}\\
&&\hskip -10mm F=\frac{t^4+u^4+v^4+w^4 -2 t^2 w^2-2 u^2 v^2}
{(t^2+u^2-v^2-w^2)( t^2-u^2+v^2-w^2)} ,
\label{3e27}\\
&&\hskip -10mm C=\frac{t^2 u^2+t^2 v^2+u^2 w^2 +v^2 w^2
- 4 t u v w}{2(t u-v w)( t v-u w)} ,
\label{3e28}\\
&&\hskip -10mm D=\frac{(t^2-w^2)^2 (u^2-v^2)^2 
\sqrt{(t^2+u^2+v^2+w^2)( t^2-u^2-v^2+w^2)(t^2 w^2-u^2 v^2)}}
{(t^2+u^2-v^2-w^2)( t^2-u^2+v^2-w^2)(t^2 u^2-v^2 w^2)(t^2 v^2-u^2 w^2)} ,
\label{3e29}
\end{eqnarray}
 where 
$C^2-D^2+E^2+F^2-2 C E F=1$ is satisfied.

If we use $P_3(0,0)=a$, $S_3(0,0)=b$, $P_0(0,0)=c$, $S_0(0,0)=d$, we have the relations 
\begin{eqnarray}
&&a^2=t^2+u^2+v^2+w^2, \quad
%\nonumber\\
%&&
b^2=2(t w +u v), \quad
%\nonumber\\
%&&
c^2=t^2-u^2-v^2+w^2,
\nonumber\\
&&d^2=2(t w -u v),
\nonumber
\end{eqnarray}
by using the addition formula. Thus, we obtain the Kummer surface relation
in the following Hudson's standard form~\cite{Hudson}
\begin{eqnarray}
&&\hskip -20mm P_1^4+S_1^4+P_2^4+S_2^4+ A_1(P_1^2 S_2^2+P_2^2 S_1^2)
+B_1(P_1^2 P_2^2+S_1^2 S_2^2)+C_1(P_1^2 S_1^2+P_2^2 S_2^2) 
\nonumber\\
&&\hskip -20mm +2D_1 P_1 P_2 S_1 S_2 =0, 
\label{3e30}\\
&&\hskip -20mm  
A_1=\frac{b^4+c^4-a^4-d^4}{a^2 d^2-b^2 c^2}, \quad B_1=\frac{a^4+c^4-b^4-d^4}{b^2 d^2-a^2 c^2}, \quad
C_1=\frac{a^4+b^4-c^4-d^4}{c^2 d^2-a^2 b^2}, 
\nonumber\\
&&\hskip -20mm  
D_1= \frac{a b c d (a^2+d^2-b^2-c^2) (b^2+d^2-a^2-c^2) (c^2+d^2-a^2-b^2) (a^2+b^2+c^2+d^2)}
{(a^2 d^2-b^2 c^2)(b^2 d^2-a^2 c^2)(c^2 d^2-a^2 b^2)}, 
\label{3e31}
\end{eqnarray}
 where 
$A_1^2+B_1^2+C_1^2-D_1^2-A_1 B_1 C_1=4$ is satisfied.

%%%%%%%%%%%%%%%%%%%%%%%%%%%%%%%%%%%%%%%%%%%%%%%%%%%%%%%%%
%%%%%%%%%%%%%%%%%%%% Section 4 %%%%%%%%%%%%%%%%%%%%%%%%%%%%%%%%%%%%%%
%%%%%%%%%%%%%%%%%%%%%%%%%%%%%%%%%%%%%%%%%%%%%%%%%%%%%%%%%%%%%%%%%
\section{The differential equation (Step I)}
\setcounter{equation}{0}
\subsection{The derivative formula}
In this section, by using  $p=S_1/P_1$, $q=S_2/P_2$ and $s=P_1/P_2$ of the form 
\begin{eqnarray}
p=\frac{S_1}{P_1}=\frac{\vartheta[\begin{array}{cc} 1 \ 1 \\ 0\ \ 1 \\ \end{array}](u,v)}
{\vartheta[\begin{array}{cc} 0 \ 0 \\ 0\ \ 1 \\ \end{array}](u,v)} , \quad
q=\frac{S_2}{P_2}=\frac{\vartheta[\begin{array}{cc} 1 \ 1 \\ 1\ \ 0 \\ \end{array}](u,v)}
{\vartheta[\begin{array}{cc} 0 \ 0 \\ 1\ \  0 \\ \end{array}](u,v)} , \quad
s=\frac{P_1}{P_2}=\frac{\vartheta[\begin{array}{cc} 0 \ 0 \\ 0\ \ 1 \\ \end{array}](u,v)}
{\vartheta[\begin{array}{cc} 0 \ 0 \\ 1\ \ 0 \\ \end{array}](u,v)}  ,
\label{4e1}
\end{eqnarray}
we first derive the differential equation of $p=S_1/P_1$ and $q=S_2/P_2$.
By noticing the relation 
\begin{eqnarray}
\dd p=\frac{1}{P_1^2}(\dd S_1 P_1 -\dd P_1 S_1),\quad  \dd q =\frac{1}{P_2^2}(\dd S_2 P_2 -\dd P_2 S_2) ,
\label{4e2}
\end{eqnarray}
we consider the combination 
\begin{eqnarray}
&&\hskip -20mm I=(\dd S_1 P_1 -\dd P_1 S_1)=\dd \vartheta[\begin{array}{cc} 1 \ 1 \\ 0\ \ 1 \\ \end{array}](u,v)
\vartheta[\begin{array}{cc} 0 \ 0 \\ 0\ \ 1 \\ \end{array}](u,v)
-\dd \vartheta[\begin{array}{cc} 0 \ 0 \\ 0\ \ 1 \\ \end{array}](u,v) 
\vartheta[\begin{array}{cc} 1 \ 1 \\ 0\ \ 1 \\ \end{array}](u,v) , 
\label{4e3}\\
&&\hskip -20mm J=(\dd S_2 P_2 -\dd P_2 S_2)=\dd \vartheta[\begin{array}{cc} 1 \ 1 \\ 1\ \ 0 \\ \end{array}](u,v)
\vartheta[\begin{array}{cc} 0 \ 0 \\ 1\ \ 0 \\ \end{array}](u,v)
-\dd \vartheta[\begin{array}{cc} 0 \ 0 \\ 1\ \ 0 \\ \end{array}](u,v) 
\vartheta[\begin{array}{cc} 1 \ 1 \\ 1\ \ 0 \\ \end{array}](u,v) .
\label{4e4}
\end{eqnarray}
We derive the derivative formula of the theta function 
 by using the addition formula of the theta function. Thus, we consider 
\begin{eqnarray}
&&\hskip -20mm I_1=\vartheta[\begin{array}{cc} 1 \ 1 \\ 0\ \ 1 \\ \end{array}](u_1+u_2,v_1+v_2)
\vartheta[\begin{array}{cc} 0 \ 0 \\ 0\ \ 1 \\ \end{array}](u_1-u_2,v_1-v_2)
\nonumber\\
&&-\vartheta[\begin{array}{cc} 0 \ 0 \\ 0\ \ 1 \\ \end{array}](u_1+u_2,v_1+v_2) 
\vartheta[\begin{array}{cc} 1 \ 1 \\ 0\ \ 1 \\ \end{array}](u_1-u_2,v_1-v_2)
\nonumber\\
&&\hskip -20mm 
=i\left(\varTheta[\begin{array}{cc} \ha \ \ha \\ 0\ \ 0 \\ \end{array}](2u_1,2v_1)
+\varTheta[\begin{array}{cc} \mha \ \mha \\ 0\ \ 0 \\ \end{array}](2u_1,2v_1)\right)
\left(\varTheta[\begin{array}{cc} \ha \ \ha \\ 0\ \ 0 \\ \end{array}](2u_2,2v_2)
-\varTheta[\begin{array}{cc} \mha \ \mha \\ 0\ \ 0 \\ \end{array}](2u_2,2v_2)\right)
\nonumber\\
&&\hskip -20mm -i\left(\varTheta[\begin{array}{cc} \ha \ \mha \\ 0\ \ 0 \\ \end{array}](2u_1,2v_1
+\varTheta[\begin{array}{cc} \mha \ \ha \\ 0\ \ 0 \\ \end{array}](2u_1,2v_1)\right)
\left(\varTheta[\begin{array}{cc} \ha \ \mha \\ 0\ \ 0 \\ \end{array}](2u_2,2v_2)
-\varTheta[\begin{array}{cc} \mha \ \ha \\ 0\ \ 0 \\ \end{array}](2u_2,2v_2)\right)  .
\nonumber
\end{eqnarray}
If we notice that
\begin{eqnarray}
&&\left(\varTheta[\begin{array}{cc} \ha \ \ha \\ 0\ \ 0 \\ \end{array}](2u_2,2v_2)
-\varTheta[\begin{array}{cc} \mha \ \mha \\ 0\ \ 0 \\ \end{array}](2u_2,2v_2)\right), 
\nonumber\\
&&\left(\varTheta[\begin{array}{cc} \ha \ \mha \\ 0\ \ 0 \\ \end{array}](2u_2,2v_2)
-\varTheta[\begin{array}{cc} \mha \ \ha \\ 0\ \ 0 \\ \end{array}](2u_2,2v_2)\right), 
\nonumber
\end{eqnarray}
are odd functions, these odd functions become zero at $u_2=0, v_2=0$.
Hence, by putting $u_2=\dd u_1, v_2=0$, and replacing 
$u_1 \rightarrow u$, $v_1 \rightarrow v$, 
we obtain  $\partial \vartheta/\partial u $ in the form
\begin{eqnarray}
&& 
\left(\frac{\partial}{\partial u} \vartheta[\begin{array}{cc} 1 \ 1 \\ 0\ \ 1 \\ \end{array}](u,v)
\vartheta[\begin{array}{cc} 0 \ 0 \\ 0\ \ 1 \\ \end{array}](u,v)
-\frac{\partial}{\partial u} \vartheta[\begin{array}{cc} 0 \ 0 \\ 0\ \ 1 \\ \end{array}](u,v) 
\vartheta[\begin{array}{cc} 1 \ 1 \\ 0\ \ 1 \\ \end{array}](u,v) \right)
\nonumber\\
&&=i \alpha \left(\varTheta[\begin{array}{cc} \ha \ \ha \\ 0\ \ 0 \\ \end{array}](2u,2v)
+\varTheta[\begin{array}{cc} \mha \ \mha \\ 0\ \ 0 \\ \end{array}](2u,2v)\right)
\nonumber\\
&&-i \beta \left(\varTheta[\begin{array}{cc} \ha \ \mha \\ 0\ \ 0 \\ \end{array}](2u,2v)
+\varTheta[\begin{array}{cc} \mha \ \ha \\ 0\ \ 0 \\ \end{array}](2u,2v)\right)
\nonumber\\
&&=i \alpha\ f(2 u, 2 v) -i \beta\ g(2 u, 2 v) ,
\label{4e5}
\end{eqnarray}
where
\begin{eqnarray}
&& f(2 u, 2 v)=\left(\varTheta[\begin{array}{cc} \ha \ \ha \\ 0\ \ 0 \\ \end{array}](2u,2v)
+\varTheta[\begin{array}{cc} \mha \ \mha \\ 0\ \ 0 \\ \end{array}](2u,2v)\right)  ,
\label{4e6}\\
&&g(2 u, 2 v)=\left(\varTheta[\begin{array}{cc} \ha \ \mha \\ 0\ \ 0 \\ \end{array}](2u,2v)
+\varTheta[\begin{array}{cc} \mha \ \ha \\ 0\ \ 0 \\ \end{array}](2u,2v)\right)  ,
\label{4e7}\\
&&\alpha=\frac{\dd}{\dd u}
\left(\varTheta[\begin{array}{cc} \ha \ \ha \\ 0\ \ 0 \\ \end{array}](u,0)
-\varTheta[\begin{array}{cc} \mha \ \mha \\ 0\ \ 0 \\ \end{array}](u, 0)\right)\biggr|_{u=0}  ,
\label{4e8}\\
&&\beta=\frac{\dd}{\dd u}
\left(\varTheta[\begin{array}{cc} \ha \ \mha \\ 0\ \ 0 \\ \end{array}](u,0)
-\varTheta[\begin{array}{cc} \mha \ \ha \\ 0\ \ 0 \\ \end{array}](u,0)\right)\biggr|_{u=0}  .
\label{4e9}
\end{eqnarray}
Similarly, putting $u_2=0$, $v_2=\dd v_1$,  and replacing 
$u_1 \rightarrow u$, $v_1 \rightarrow v$, we obtain  $\partial \vartheta/\partial v$ 
in the form
\begin{eqnarray}
&& 
\left(\frac{\partial}{\partial v} \vartheta[\begin{array}{cc} 1 \ 1 \\ 0\ \ 1 \\ \end{array}](u,v)
\vartheta[\begin{array}{cc} 0 \ 0 \\ 0\ \ 1 \\ \end{array}](u,v)
-\frac{\partial}{\partial v} \vartheta[\begin{array}{cc} 0 \ 0 \\ 0\ \ 1 \\ \end{array}](u,v) 
\vartheta[\begin{array}{cc} 1 \ 1 \\ 0\ \ 1 \\ \end{array}](u,v) \right)
\nonumber\\
&&=i \gamma \left(\varTheta[\begin{array}{cc} \ha \ \ha \\ 0\ \ 0 \\ \end{array}](2u,2v)
+\varTheta[\begin{array}{cc} \mha \ \mha \\ 0\ \ 0 \\ \end{array}](2u,2v)\right)
\nonumber\\
&&-i \delta \left(\varTheta[\begin{array}{cc} \ha \ \mha \\ 0\ \ 0 \\ \end{array}](2u,2v)
+\varTheta[\begin{array}{cc} \mha \ \ha \\ 0\ \ 0 \\ \end{array}](2u,2v)\right)
\nonumber\\
&&=i \gamma\ f(2 u, 2 v) -i \delta\ g(2 u, 2 v) ,
\label{4e10}\\
%\end{eqnarray}
%
%where
%
%\begin{eqnarray}
&&\gamma=\frac{\dd}{\dd v}
\left(\varTheta[\begin{array}{cc} \ha \ \ha \\ 0\ \ 0 \\ \end{array}](0,v)
-\varTheta[\begin{array}{cc} \mha \ \mha \\ 0\ \ 0 \\ \end{array}](0, v)\right)\biggr|_{v=0}  ,
\label{4e11}\\
&&\delta=\frac{\dd}{\dd v}
\left(\varTheta[\begin{array}{cc} \ha \ \mha \\ 0\ \ 0 \\ \end{array}](0,v)
-\varTheta[\begin{array}{cc} \mha \ \ha \\ 0\ \ 0 \\ \end{array}](0, v)\right)\biggr|_{v=0}  .
\label{4e12}
\end{eqnarray}
Thus, we obtain the derivative formula from the addition formula of the theta function 
in the form 
\begin{eqnarray}
&&I=(\dd S_1 P_1 -\dd P_1 S_1)=\dd u (\partial_u S_1 P_1-\partial_u P_1 S_1)
+\dd v (\partial_v S_1 P_1-\partial_v P_1 S_1)
\nonumber\\
&&=i f(2u,2v) (\alpha \dd u+\gamma \dd v)  -i g(2u,2v) (\beta \dd u+\delta \dd v)  .
\label{4e13}
\end{eqnarray}
We can obtain $J$ by the replacement of
$u \rightarrow u+1$, $v \rightarrow v+1$ in Eq.(\ref{4e13}). In the right-hand side
of Eq.(\ref{4e13}), we obtain   

\begin{eqnarray}
&&\hskip -20mm I=(\dd S_1 P_1 -\dd P_1 S_1)=\dd \vartheta[\begin{array}{cc} 1 \ 1 \\ 0\ \ 1 \\ \end{array}](u,v)
\vartheta[\begin{array}{cc} 0 \ 0 \\ 0\ \ 1 \\ \end{array}](u,v)
-\dd \vartheta[\begin{array}{cc} 0 \ 0 \\ 0\ \ 1 \\ \end{array}](u,v) 
\vartheta[\begin{array}{cc} 1 \ 1 \\ 0\ \ 1 \\ \end{array}](u,v) 
\nonumber\\
&&\hskip -20mm \rightarrow
\dd \vartheta[\begin{array}{cc} 1 \ 1 \\ 1\ \ 2 \\ \end{array}](u,v)
\vartheta[\begin{array}{cc} 0 \ 0 \\ 1\ \ 2 \\ \end{array}](u,v)
-\dd \vartheta[\begin{array}{cc} 0 \ 0 \\ 1\ \ 2 \\ \end{array}](u,v) 
\vartheta[\begin{array}{cc} 1 \ 1 \\ 1\ \ 2 \\ \end{array}](u,v)\nonumber\\
 &&\hskip -20mm =(-1)\left(\dd \vartheta[\begin{array}{cc} 1 \ 1 \\ 1\ \ 0 \\ \end{array}](u,v)
\vartheta[\begin{array}{cc} 0 \ 0 \\ 1\ \ 0 \\ \end{array}](u,v)
-\dd \vartheta[\begin{array}{cc} 0 \ 0 \\ 1\ \ 0 \\ \end{array}](u,v) 
\vartheta[\begin{array}{cc} 1 \ 1 \\ 1\ \ 0 \\ \end{array}](u,v)\right)=-J   .
\label{4e14}
\end{eqnarray}
While, in the right-hand side of Eq.(\ref{4e13}), we obtain 
\begin{eqnarray}
&&i f(2u,2v) (\alpha \dd u+\gamma \dd v)  -i g(2u,2v) (\beta \dd u+\delta \dd v)  
\nonumber\\
&&=i  \left(\varTheta[\begin{array}{cc} \ha \ \ha \\ 0\ \ 0 \\ \end{array}](2u,2v)
+\varTheta[\begin{array}{cc} \mha \ \mha \\ 0\ \ 0 \\ \end{array}](2u,2v)\right)
(\alpha \dd u+\gamma \dd v) 
\nonumber\\
&&-i \left(\varTheta[\begin{array}{cc} \ha \ \mha \\ 0\ \ 0 \\ \end{array}](2u,2v)
+\varTheta[\begin{array}{cc} \mha \ \ha \\ 0\ \ 0 \\ \end{array}](2u,2v)\right)
 (\beta \dd u+\delta \dd v) 
\nonumber\\
&&\rightarrow 
i  \left(\varTheta[\begin{array}{cc} \ha \ \ha \\ 2\ \ 2 \\ \end{array}](2u,2v)
+\varTheta[\begin{array}{cc} \mha \ \mha \\ 2\ \ 2 \\ \end{array}](2u,2v)\right)
(\alpha \dd u+\gamma \dd v) 
\nonumber\\
&&-i \left(\varTheta[\begin{array}{cc} \ha \ \mha \\ 2\ \ 2 \\ \end{array}](2u,2v)
+\varTheta[\begin{array}{cc} \mha \ \ha \\ 2\ \ 2 \\ \end{array}](2u,2v)\right)
 (\beta \dd u+\delta \dd v) 
\nonumber\\
&&=-i  \left(\varTheta[\begin{array}{cc} \ha \ \ha \\ 0\ \ 0 \\ \end{array}](2u,2v)
+\varTheta[\begin{array}{cc} \mha \ \mha \\ 0\ \ 0 \\ \end{array}](2u,2v)\right)
(\alpha \dd u+\gamma \dd v) 
\nonumber\\
&&-i \left(\varTheta[\begin{array}{cc} \ha \ \mha \\ 0\ \ 0 \\ \end{array}](2u,2v)
+\varTheta[\begin{array}{cc} \mha \ \ha \\ 0\ \ 0 \\ \end{array}](2u,2v)\right)
 (\beta \dd u+\delta \dd v) \nonumber\\
&&=-i f(2u, 2v)(\alpha \dd u+\gamma \dd v) 
-i g(2u, 2v)  (\beta \dd u+\delta \dd v) .
\label{4e15}
\end{eqnarray}
Thus, we obtain the derivative formula of $J$ in the form\footnote{In the 
G\"{o}pel's expression, there is the overall minus sign.}
\begin{eqnarray}
&&J=(\dd S_2 P_2 -\dd P_2 S_2)=\dd u (\partial_u S_2 P_2-\partial_u P_2 S_2)
+\dd v (\partial_v S_2 P_2-\partial_v P_2 S_2), 
\nonumber\\
&&=i f(2u,2v) (\alpha \dd u+\gamma \dd v)  +i g(2u,2v) (\beta \dd u+\delta \dd v) . 
\label{4e16}
\end{eqnarray}
Therefore, the necessary derivative formulae provide 
\begin{eqnarray}
&&(\dd S_1 P_1 -\dd P_1 S_1)=i f(2u,2v) (\alpha \dd u+\gamma \dd v)  -i g(2u,2v) (\beta \dd u+\delta \dd v)  ,
\label{4e17}\\
&&(\dd S_2 P_2 -\dd P_2 S_2)=i f(2u,2v) (\alpha \dd u+\gamma \dd v)  +i g(2u,2v) (\beta \dd u+\delta \dd v) .
\label{4e18} 
\end{eqnarray}
Next, we express $ f(2u,2v)$, $g(2u,2v)$ with the product of the original   
theta function 
$P_0$, $P_1$, $\cdots$, $S_2$, $S_3$. Using Eq.(\ref{3e17}), we obtain 
\begin{eqnarray}
&&\hskip -10mm P_3 S_3=\vartheta[\begin{array}{cc} 1 \ 1 \\ 0\ \ 0 \\ \end{array}](u, v) 
\vartheta[\begin{array}{cc} 0 \ 0 \\ 0\ \ 0 \\ \end{array}](u, v)
\nonumber\\
&&\hskip -10mm =\left(\varTheta[\begin{array}{cc} \ha \ \ha \\ 0\ \ 0 \\ \end{array}](2u, 2v) 
+\varTheta[\begin{array}{cc} \mha \ \mha \\ 0\ \ 0 \\ \end{array}](2u, 2v)\right)
\varTheta[\begin{array}{cc} \ha \ \ha \\ 0\ \ \ 0 \\ \end{array}](0, 0)
\nonumber\\
&&\hskip -10mm 
+\left(\varTheta[\begin{array}{cc} \ha \ \mha \\ 0\ \ 0 \\ \end{array}](2u, 2v) 
+\varTheta[\begin{array}{cc} \mha \ \ha \\ 0\ \ 0 \\ \end{array}](2u, 2v)\right)
\varTheta[\begin{array}{cc} \ha \ \mha \\ 0\ \ \ 0 \\ \end{array}](0, 0)
\nonumber\\
&&\hskip -10mm =\varTheta[\begin{array}{cc} \ha \ \ha \\ 0\ \ \ 0 \\ \end{array}](0, 0) f(2u, 2v)
+\varTheta[\begin{array}{cc} \ha \ \mha \\ 0\ \ \ 0 \\ \end{array}](0, 0) g(2u,2v)  .
\label{4e19}
\end{eqnarray} 
By replacing $u \rightarrow u+1, v \rightarrow v+1$ in Eq.(\ref{4e19}), we obtain
\begin{eqnarray}
&&\hskip -10mm P_0 S_0=\vartheta[\begin{array}{cc} 1 \ 1 \\ 1\ \ 1 \\ \end{array}](u, v) 
\vartheta[\begin{array}{cc} 0 \ 0 \\ 1\ \ 1 \\ \end{array}](u, v)
\nonumber\\
&&\hskip -10mm =\left(\varTheta[\begin{array}{cc} \ha \ \ha \\ 2\ \ 2 \\ \end{array}](2u, 2v) 
+\varTheta[\begin{array}{cc} \mha \ \mha \\ 2\ \ 2 \\ \end{array}](2u, 2v)\right)
\varTheta[\begin{array}{cc} \ha \ \ha \\ 0\ \ \ 0 \\ \end{array}](0, 0)
\nonumber\\
&&\hskip -10mm 
+\left(\varTheta[\begin{array}{cc} \ha \ \mha \\ 2\ \ 2 \\ \end{array}](2u, 2v) 
+\varTheta[\begin{array}{cc} \mha \ \ha \\ 2\ \ 2 \\ \end{array}](2u, 2v)\right)
\varTheta[\begin{array}{cc} \ha \ \mha \\ 0\ \ \ 0 \\ \end{array}](0, 0)
\nonumber\\
&&\hskip -10mm =-\varTheta[\begin{array}{cc} \ha \ \ha \\ 0\ \ \ 0 \\ \end{array}](0, 0) f(2u, 2v)
+\varTheta[\begin{array}{cc} \ha \ \mha \\ 0\ \ \ 0 \\ \end{array}](0, 0) g(2u,2v)  .
\label{4e20}
\end{eqnarray} 
Thus, we can express $f(2u, 2v)$ and $g(2u, 2v)$ with the product of the original 
theta functions in the form
\begin{eqnarray}
&&f(2u,2v)=\frac{P_3 S_3-P_0 S_0}{2 \varTheta[\begin{array}{cc} \ha \ \ha \\ 0\ \ \ 0 \\ \end{array}](0, 0)}, \quad
g(2u,2v)=\frac{P_3 S_3+P_0 S_0}{2 \varTheta[\begin{array}{cc} \ha \ \mha \\ 0\ \ \ 0 \\ \end{array}](0, 0)}  .
\label{4e21}
\end{eqnarray}
%
%%%%%%%%%%%%%%%%%%%%%%%%%%%%%%%%%%%
\subsection{The differential equation (Step I)}
By using Eqs.(\ref{4e17}), (\ref{4e18}), and (\ref{4e21}), we obtain
the following differential equation 
\begin{eqnarray}
&&\hskip -20mm (\dd S_1 P_1 -\dd P_1 S_1)
=i\frac{P_3 S_3-P_0 S_0}{2 \varTheta[\begin{array}{cc} \ha \ \ha \\ 0\ \ \ 0 \\ \end{array}](0, 0)}
(\alpha \dd u+\gamma \dd v)  
-i\frac{P_3 S_3+P_0 S_0}{2 \varTheta[\begin{array}{cc} \ha \ \mha \\ 0\ \ \ 0 \\ \end{array}](0, 0)}
 (\beta \dd u+\delta \dd v) , 
\label{4e22}\\
&&\hskip -20mm (\dd S_2 P_2 -\dd P_2 S_2)
=i \frac{P_3 S_3-P_0 S_0}{2 \varTheta[\begin{array}{cc} \ha \ \ha \\ 0\ \ \ 0 \\ \end{array}](0, 0)}
 (\alpha \dd u+\gamma \dd v)  
+i \frac{P_3 S_3+P_0 S_0}{2 \varTheta[\begin{array}{cc} \ha \ \mha \\ 0\ \ \ 0 \\ \end{array}](0, 0)}
(\beta \dd u+\delta \dd v) .
\label{4e23} 
\end{eqnarray}
We can verify that Eq.(\ref{4e23}) is obtained from Eq.(\ref{4e22}) by replacing
$u \rightarrow u+1, v\rightarrow v+1$.\footnote{In the G\"{o}pel's expression, in the right-hand side of Eq.(\ref{4e23}), there is the overall minus sign.}
Under such replacement, $(\dd S_1 P_1 -\dd P_1 S_1) \rightarrow -(\dd S_2 P_2 -\dd P_2 S_2)$,
$(P_3 S_3+P_0 S_0) \rightarrow (P_3 S_3+P_0 S_0) $, $(P_3 S_3-P_0 S_0) \rightarrow -(P_3 S_3-P_0 S_0) $, which provides the proof that Eq.(\ref{4e23}) is obtained from Eq.(\ref{4e22}).   
Thus, we obtain the following differential equation 
\begin{eqnarray}
&&\frac{(\dd S_1 P_1 -\dd P_1 S_1)-(\dd S_2 P_2 -\dd P_2 S_2)}{P_3 S_3+P_0 S_0}=
\frac{-i(\beta \dd u+\delta \dd v)}
{\varTheta[\begin{array}{cc} \ha \ \mha \\ 0\ \ \ 0 \\ \end{array}](0, 0)}
=\dd \mu ,
\label{4e24}\\
&&\frac{(\dd S_1 P_1 -\dd P_1 S_1)+(\dd S_2 P_2 -\dd P_2 S_2)}{P_3 S_3-P_0 S_0}=
\frac{i(\alpha \dd u+\gamma \dd v)}
{\varTheta[\begin{array}{cc} \ha \ \ha \\ 0\ \ \ 0 \\ \end{array}](0, 0)}
=\dd \nu  .
\label{4e25}
\end{eqnarray}
Which can be expressed in the form
\begin{eqnarray}
&&\frac{(P_1^2 \dd p- P_2^2 \dd q )/P_1 P_2}{(P_3 S_3+P_0 S_0)/P_1 P_2}=
\frac{s \dd p -\dd q/s}{\varphi}=\dd \mu,
\label{4e26}\\
&&\frac{(P_1^2 \dd p +P_2^2 \dd q )/P_1 P_2}{(P_3 S_3-P_0 S_0)/P_1 P_2}=
\frac{s \dd p +\dd q/s}{\psi}=\dd \nu  ,
\label{4e27}\\
&&s=\frac{P_1}{P_2}, \quad \psi=\frac{P_3 S_3-P_0 S_0}{P_1 P_2}, \quad \varphi=\frac{P_3 S_3+P_0 S_0}{P_1 P_2}.
\label{4e28}
\end{eqnarray}
Therefore, we obtain the starting differential equation\footnote{In G\"{o}pel's 
expression, $(\pm 1)$ sign of $s \dd p \mp \dd q/s$ is in the opposite}
\begin{eqnarray}
&&\frac{s \dd p -\dd q/s}{\varphi}=\dd \mu, \quad
\frac{s \dd p +\dd q/s}{\psi}=\dd \nu,
\label{4e29}\\
&&p=\frac{S_1}{P_1}, \quad q=\frac{S_2}{P_2}, \quad s=\frac{P_1}{P_2}
\nonumber\\
&&\varphi=\frac{P_3 S_3+P_0 S_0}{P_1 P_2}, \quad
\psi=\frac{P_3 S_3-P_0 S_0}{P_1 P_2}.
\label{4e30}
\end{eqnarray}
%

%%%%%%%%%%%%%%%%%%%%%%%%%%%%%%%%%%%%%%%%%%%%%%%%%%%%%%%%%%%%%%%%%%%%%%%%%%%%%%
%%%%%%%%%%%%%%%%%%%%%%%%%%   Section 5  %%%%%%%%%%%%%%%%%%%%%%%%%%%%%%%%%%%%%%%%%
%%%%%%%%%%%%%%%%%%%%%%%%%%%%%%%%%%%%%%%%%%%%%%%%%%%%%%%%%%%%%%%%%%%%%%%%%%%%%
\section{The differential equation with $p$ and $q$ (Step II)}
\setcounter{equation}{0}
\subsection{The elimination of $s$-dependence in 
$\displaystyle{s \dd p \mp \frac{\dd q}{s}}$}
By dividing $P_1^2 P_2^2$ in the  Kummer surface relation Eq.(\ref{3e25}), we obtain 
\begin{eqnarray}
(1-2E p^2+p^4)s^4+(1-2Eq^2+q^4)-2\Big(F(1+p^2 q^2)-C(p^2+q^2)+2Dpq\Big)s^2=0  ,
\label{5e1}
\end{eqnarray}
which provides
\begin{eqnarray}
(1-2E p^2+p^4)s^2+\frac{(1-2Eq^2+q^4)}{s^2}=2\Big(F(1+p^2 q^2)-C(p^2+q^2)+2Dpq\Big)  .
\label{5e2}
\end{eqnarray}
This is given in the form
\begin{eqnarray}
&&\left(\sqrt{1-2E p^2+p^4}\ s \pm \frac{\sqrt{1-2Eq^2+q^4}}{s}\right)^2
\nonumber\\
&&=2\Big(F(1+p^2 q^2)-C(p^2+q^2)+2Dpq\Big) \pm 2 \sqrt{1-2E p^2+p^4} \sqrt{1-2Eq^2+q^4} .
\label{5e3}
\end{eqnarray}
Introducing the function $\Delta$ of the form 
$\Delta(x)=\sqrt{1- 2 E x^2 +x^4}$, we can express the above in the form
\begin{eqnarray}
&&\left(\Delta(p)\ s \pm \frac{\Delta(q)}{s}\right)^2=2G(p,q) \pm 2\Delta(p) \Delta(q),
\label{5e4}\\
&& G(p,q)=\Big(F(1+p^2 q^2)-C(p^2+q^2)+2Dpq\Big) .
\label{5e5}
\end{eqnarray}
Thus, we obtain the relation
\begin{eqnarray}
&&\Delta(p)\ s \pm \frac{\Delta(q)}{s}=\sqrt{2G(p,q) \pm 2\Delta(p) \Delta(q)}  .
\label{5e6}
\end{eqnarray}
Therefore, in the combination of $(\Delta(p) s \pm \Delta(q)/s)$, $s$-dependence 
is eliminated.
In order to eliminate $s$-dependence of the differential, we rearrange 
$(\Delta(p) s \pm \Delta(q)/s)$ in the form 
\begin{eqnarray}
&&s\ \dd p -\frac{\dd q}{s}
\nonumber\\
&&
=\left(s \Delta(p)-\frac{\Delta(q)}{s}\right)
\frac{1}{2} \left(\frac{\dd p }{\Delta(p)} +\frac{\dd q }{\Delta(q) }\right)
+\left(s \Delta(p)+\frac{\Delta(q)}{s}\right)
\frac{1}{2} \left(\frac{\dd p }{\Delta(p)} -\frac{\dd q }{\Delta(q) }\right)
\nonumber\\
&&
=\sqrt{2G(p,q) -2\Delta(p) \Delta(q)}\frac{1}{2} \left(\frac{\dd p }{\Delta(p)} +\frac{\dd q }{\Delta(q) }\right)
\nonumber\\
&&
+\sqrt{2G(p,q) +2\Delta(p) \Delta(q)}
\frac{1}{2} \left(\frac{\dd p }{\Delta(p)} -\frac{\dd q }{\Delta(q) }\right) , 
\label{5e7}\\
&&s\ \dd p +\frac{\dd q}{s}
\nonumber\\
&&
=\left(s \Delta(p)+\frac{\Delta(q)}{s}\right)
\frac{1}{2} \left(\frac{\dd p }{\Delta(p)} +\frac{\dd q }{\Delta(q) }\right)
+\left(s \Delta(p)-\frac{\Delta(q)}{s}\right)
\frac{1}{2} \left(\frac{\dd p }{\Delta(p)} -\frac{\dd q }{\Delta(q) }\right)
\nonumber\\
&&
=\sqrt{2G(p,q) +2\Delta(p) \Delta(q)}
\frac{1}{2} \left(\frac{\dd p }{\Delta(p)} +\frac{\dd q }{\Delta(q) }\right)
\nonumber\\
&& +\sqrt{2G(p,q)-2\Delta(p) \Delta(q)}
\frac{1}{2} \left(\frac{\dd p }{\Delta(p)} -\frac{\dd q }{\Delta(q) }\right) .
\label{5e8}
\end{eqnarray}
Thus, we can express $\displaystyle{s \dd p \mp \frac{\dd q}{s}}$ only with $p$ and $q$.

%%%%%%%%%%%%%%%%%%%%%%%%%%%
\subsection{The expression of  $\psi$ and $\varphi$ with $p$ and $q$ } 

\noindent
{\bf a) The expression of $(P_3^2 S_3^2+P_0^2 S_0^2)/P_1^2 P_2^2$ with $p$ and $q$}\\
Next, we express $(P_3^2 S_3^2+P_0^2 S_0^2)/P_1^2 P_2^2$ as the 
function of $p$ and $q$.
By using the addition formula, we first express 
$P_3^2, S_3^2, P_0^2, S_0^2, P_1^2, P_2^2$, as the function of $T, U, V, W$ via Eqs.(\ref{2e5})-(\ref{2e8}); in addition, 
by using Eq.(\ref{3e3}) and Eq.(\ref{3e4}), $P_3^3 S_3^2$, $P_0^3 S_0^2$
and $P_1^3 S_1^2$ are expressed as the function of $P_1^2, S_1^2, P_2^2, S_2^2$.
Thus, $(P_3^2 S_3^2+P_0^2 S_0^2)/P_1^2 P_2^2$ is expressed
as the function of $p^2, q^2, s^2$ in the form
\begin{eqnarray}
&&\hskip -20mm \frac{1}{2}(\varphi^2+\psi^2)=\frac{P_3^2 S_3^2+P_0^2 S_0^2}{P_1^2 P_2^2} 
=\frac{(t^2 u^2-v^2 w^2)(t^2 v^2 -u^2 w^2)}{(t^2-w^2)^2 (u^2-v^2)^2}
\nonumber\\
&&\hskip -20mm \times \left(2C\left((1-2 E p^2+p^4)s^2 +\frac{(1-2 E q^2+q^4)}{s^2}\right)
-4E(1+p^2q^2)+4 (p^2+q^2)\right) .
\label{5e9}
\end{eqnarray}
The $s$-dependence is eliminated by using the Kummer's quartic relation Eq.(\ref{5e2})
of the form
\begin{eqnarray}
(1-2E p^2+p^4)s^2+\frac{(1-2Eq^2+q^4)}{s^2}=2\Big(F(1+p^2 q^2)-C(p^2+q^2)+2Dpq\Big)=2G(p,q)  .
\nonumber
\end{eqnarray}
Thus, we obtain the $s$-independent expression of $(P_3^2 S_3^2+P_0^2 S_0^2)/P_1^2 P_2^2$ 
in the form
\begin{eqnarray}
&&\frac{1}{2}(\varphi^2+\psi^2)=\frac{P_3^2 S_3^2+P_0^2 S_0^2}{P_1^2 P_2^2} 
=\frac{(t^2 u^2-v^2 w^2)(t^2 v^2 -u^2 w^2)}{(t^2-w^2)^2 (u^2-v^2)^2}
\nonumber\\
&&\times \left(4C(F(1+p^2 q^2)-C(p^2+q^2)+2 Dpq)-4E(1+p^2q^2)+4 (p^2+q^2)\right)
\nonumber\\
&&= b (1+p^2 q^2) -a (p^2+q^2) +2 c pq  ,
\label{5e10}
\end{eqnarray}
where 
\begin{eqnarray}
&&\hskip -20mm a=\frac{C^2-1}{\sqrt{(C^2-1)(E^2-1)}},\ b=\frac{CF-E}{\sqrt{(C^2-1)(E^2-1)}}, \ 
 c=\frac{CD}{\sqrt{(C^2-1)(E^2-1)}}   ,
\nonumber
\end{eqnarray}
and we used
$\displaystyle{ \sqrt{(C^2-1)(E^2-1)}
=\frac{(t^2-w^2)^2 (u^2-v^2)^2}{4(t^2 u^2-v^2 w^2)(t^2 v^2-u^2 w^2)}}$.
Thus, we obtain
\begin{eqnarray}
&&\frac{1}{2}(\varphi^2+\psi^2)=\frac{P_3^2 S_3^2+P_0^2 S_0^2}{P_1^2 P_2^2} = b (1+p^2 q^2) -a (p^2+q^2) +2 c pq .
\label{5e11}
\end{eqnarray}
%
%%%%%%%%%%%%%%%%%%%%%%%
\noindent 
{\bf b) The expression of $P_3^2 S_3^2P_0^2 S_0^2/P_1^4 P_2^4$ with $p$ and $q$}\\
Similarly, we calculate 
\begin{eqnarray}
&&\varphi \psi=\frac{P_3^2 S_3^2-P_0^2 S_0^2}{P_1^2 P_2^2} 
=\frac{(tu+v w)(t v+uw)}{(t^2-w^2) (u^2-v^2)}
\nonumber\\
&&\times \left( (1-2 E p^2+p^4)s^2 -\frac{(1-2 E q^2+q^4)}{s^2})\right)
\nonumber\\
&&=\frac{1}{2\sqrt{E^2-1}} \left( (1-2 E p^2+p^4)s^2 -\frac{(1-2 E q^2+q^4)}{s^2})\right)  ,
\label{5e12}
\end{eqnarray}
where
$\displaystyle{2\sqrt{E^2-1} =\frac{(t^2-w^2) (u^2-v^2)}{(tu+v w)(t v+uw)}}$ is used. 
Using Eq.(\ref{5e2}) and Eq.(\ref{5e12}), we obtain
\begin{eqnarray}
&&\left( (1-2 E p^2+p^4)s^2 +\frac{(1-2 E q^2+q^4)}{s^2})\right)=2 G(p,q)  ,
\label{5e13}\\
&&\left( (1-2 E p^2+p^4)s^2 -\frac{(1-2 E q^2+q^4)}{s^2})\right)
=2\sqrt{E^2-1}\varphi \psi  , 
\label{5e14}
\end{eqnarray}
which provides 
\begin{eqnarray}
&&s^2=\frac{G(p,q)+\sqrt{E^2-1}\varphi \psi}{1-2 E p^2+p^4}, \quad
\frac{1}{s^2}=\frac{G(p,q)-\sqrt{E^2-1}\varphi \psi}{1-2 E q^2+q^4}  ,
\label{5e15}
\end{eqnarray}
Multiplying the first and the second term, we obtain 
\begin{eqnarray}
&&1=\frac{G(p,q)^2-(E^2-1)\varphi^2 \psi^2}{(1-2 E p^2+p^4)(1-2 Eq^2+q^4)}  .
\nonumber
\end{eqnarray}
Thus, we obtain $\varphi^2 \psi^2$ expressed with $p$ and $q$ 
\begin{eqnarray}
&&(E^2-1) \varphi^2 \psi^2=G(p,q)^2-(1-2 E p^2+p^4)(1-2 Eq^2+q^4)
=G(p,q)^2-\Delta(p)^2 \Delta(q)^2.
\nonumber
\end{eqnarray}
Hence, we obtain $\varphi \psi$ as the function of $p$ and $q$ 
\begin{eqnarray}
\varphi \psi=\frac{P_3^2 S_3^2-P_0^2 S_0^2}{P_1^2 P_2^2} 
=\sqrt{\frac{G(p,q)^2-\Delta(p^2 \Delta(q)^2}{E^2-1}}   .
\label{5e16}
\end{eqnarray}
Using Eq.(\ref{5e11}) and Eq.(\ref{5e16}), we obtain
\begin{eqnarray}
&&K_1=\frac{P_3^2 S_3^2+P_0^2 S_0^2}{P_1^2 P_2^2} =
\frac{1}{2}(\varphi^2+\psi^2)= b (1+p^2 q^2) -a (p^2+q^2) +2c pq ,
\label{5e17}\\
&&K_2=\frac{P_3^2 S_3^2-P_0^2 S_0^2}{P_1^2 P_2^2} =\varphi \psi
=\sqrt{\frac{G(p,q)^2-\Delta(p^2 \Delta(q)^2}{E^2-1}}  .
\label{5e18}
\end{eqnarray}
Thus, we obtain $P_3^2 S_3^2P_0^2 S_0^2/P_1^4 P_2^4$ in the form
\begin{eqnarray}
&&\frac{P_3^2 S_3^2P_0^2 S_0^2}{P_1^4 P_2^4}=\frac{K_1^2-K_2^2}{4}
\nonumber\\
&&=\frac{1}{4} \left( \left(b (1+p^2 q^2) -a (p^2+q^2) +2c pq \right)^2
-\frac{ G(p,q)^2-\Delta(p)^2 \Delta(q)^2}{E^2-1} \right)
\nonumber\\
&&=\left(\frac{1}{2}b_1(1+p^2 q^2)-c_1 p q\right)^2  ,
\label{5e19}\\
&&b_1=\frac{D}{\sqrt{(C^2-1)(E^2-1)}}, \quad c_1
=\frac{CE-F}{\sqrt{(C^2-1)(E^2-1)}}  ,
\nonumber
\end{eqnarray}
where we used $C^2-D^2+E^2+F^2-2 C E F=1$.
Thus,  we obtain
\begin{eqnarray}
\frac{1}{4}(\varphi^2-\psi^2)=\frac{P_3 S_3 P_0 S_0}{P_1^2 P_2^2}
=\pm \left(\frac{1}{2}b_1(1+p^2 q^2)-c_1 p q\right) .
\label{5e20}
\end{eqnarray}
For the sign ambiguity, we take $(+1)$ sign.
Therefore, we obtain $\varphi^2$ and $\psi^2$ as the function of $p$ and $q$
\begin{eqnarray}
&&\hskip -15mm \varphi^2=\frac{P_3^2 S_3^2+P_0^2 S_0^2+ 2P_3 S_3 P_0 S_0 }{P_1^2 P_2^2} 
=(b+b_1)(1+p^2 q^2)-a (p^2+q^2)+2(c-c_1)pq  ,
\label{5e21}\\
&&\hskip -15mm \psi^2=\frac{P_3^2 S_3^2+P_0^2 S_0^2- 2P_3 S_3 P_0 S_0 }{P_1^2 P_2^2} 
=(b-b_1)(1+p^2 q^2)-a (p^2+q^2)+2(c+c_1)pq  .
\label{5e22}
\end{eqnarray}
Thus, $\varphi$ and $\psi$ are expressed by the function of $p$ and $q$ in the form
\begin{eqnarray}
&&
%\hskip-30mm 
\varphi=\sqrt{(b+b_1)(1+p^2 q^2)-a(p^2+q^2)+2 (c-c_1)p q}, 
\label{5e23}\\
&&
%\hskip-30mm 
\psi=\sqrt{(b-b_1)(1+p^2 q^2)-a(p^2+q^2)+2 (c+c_1)p q} . 
\label{5e24}
\end{eqnarray}
%
%%%%%%%%%%%%%%%%%%%%%%%%%%
\subsection{The differential equation with $p$ and $q$ (Step II)}
Combining the previous results, we obtain the following differential equation, expressed 
with only $p$ and $q$ 
\begin{eqnarray}
&& \frac{1}{\varphi} \left\{ \sqrt{2 G(p,q)-2 \triangle(p) \triangle(q)} 
\frac{1}{2}\left( \frac{\dd p}{\triangle(p)}+\frac{\dd q}{\triangle(q)}\right) \right.
\nonumber\\
&& \left. +\sqrt{2 G(p,q)+2 \triangle(p) \triangle(q)} 
\frac{1}{2}\left( \frac{\dd p}{\triangle(p)}-\frac{\dd q}{\triangle(q)}\right)\right\}
=\dd \mu, 
\label{5e25}\\
&& \frac{1}{\psi} \left\{ \sqrt{2 G(p,q)+2 \triangle(p) \triangle(q)} 
\frac{1}{2}\left( \frac{\dd p}{\triangle(p)}+\frac{\dd q}{\triangle(q)}\right) \right.
\nonumber\\
&& \left. +\sqrt{2 G(p,q)-2 \triangle(p) \triangle(q)} 
\frac{1}{2}\left( \frac{\dd p}{\triangle(p)}-\frac{\dd q}{\triangle(q)}\right)\right\}
=\dd \nu, 
\label{5e26}\\
&&
\varphi=\sqrt{(b+b_1)(1+p^2 q^2)-a(p^2+q^2)+2 (c-c_1)p q}, 
\label{5e27}\\
&&
\psi=\sqrt{(b-b_1)(1+p^2 q^2)-a(p^2+q^2)+2 (c+c_1)p q}, 
\label{5e28}\\
&&
G(p,q)=F(1+p^2 q^2)-C(p^2+q^2)+2D pq, 
\quad \triangle(x)=\sqrt{1-2 E x^2+ x^4}  ,
\label{5e29}\\
&&
b=\frac{C F -E}{ \sqrt{(C^2-1)(E^2-1)}}, \ 
a=\frac{C^2-1}{ \sqrt{(C^2-1)(E^2-1)}}, \ 
c=\frac{CD}{ \sqrt{(C^2-1)(E^2-1)}}, 
\nonumber\\
&&
b_1=\frac{D}{ \sqrt{(C^2-1)(E^2-1)}}, \ 
c_1=\frac{C E -F}{ \sqrt{(C^2-1)(E^2-1)}}  .
\nonumber
\end{eqnarray}
%

 %%%%%%%%%%%%%%%%%%%%%%%%%%%%%%%%%%%%%%%%%%%%%%%%%%%%%%%%%%%%%%%%%%%%%%%%%%%%%%
%%%%%%%%%%%%%%%%%%%%%%%%%%   Section 6  %%%%%%%%%%%%%%%%%%%%%%%%%%%%%%%%%%%%%%%%%
%%%%%%%%%%%%%%%%%%%%%%%%%%%%%%%%%%%%%%%%%%%%%%%%%%%%%%%%%%%%%%%%%%%%%%%%%%%%%
%%%%%%%%%%%%%%%%%%%%%%%%%%%%%%%%%%%%%%%%%%%%%%%%%%%%%%%%%%%%
\section{The differential equation with $y$ and $z$ (Step III)}
\setcounter{equation}{0}
%%%%%%%%%%%%%%%%%%%%%%%%%%%%%%%
\subsection{The change of  the functions from $p$, $q$ to $y$, $z$}
Next, we change the function from $p$, $q$ into $y$, $z$ in such a way as the 
differential equation becomes separable. 
Because the combination of 
$\displaystyle{\frac{\dd p}{\triangle(p)} \pm \frac{\dd q}{\triangle(q)}}$ 
emerges, we change from $p$, $q$ to $y$, $z$ in such a way 
as $y$, $z$ satisfy the following differential equation
\begin{eqnarray}
\frac{\dd p}{\triangle(p)}=\frac{\dd y}{\triangle(y)}+ \frac{\dd z}{\triangle(z)}, \quad
\frac{\dd q}{\triangle(q)}=\frac{\dd y}{\triangle(y)}- \frac{\dd z}{\triangle(z)} .
\label{6e1}
\end{eqnarray}
This is the differential equation, which provides the addition 
formula of the elliptic function for
the Jacobi type the elliptic curve $y^2=1-2 E x^2+x^4=\Delta(x)^2$.
Thus, the addition formula gives 
\begin{eqnarray}
p=\frac{y \Delta(z)+z \Delta(y)}{1-y^2 z^2}, \quad q=\frac{y \Delta(z)-z \Delta(y)}{1-y^2 z^2} . 
\label{6e2}
\end{eqnarray}
For more general Jacobi type elliptic curve $y^2=1+\lambda_2 x^2+\lambda_4 x^4$,
we put $\Delta(x)=\sqrt{1+\lambda_2 x^2+\lambda_4 x^4}$ and the 
differential equation 
$\displaystyle{\frac{\dd p}{\triangle(p)}=\frac{\dd y}{\triangle(y)}+ \frac{\dd z}{\triangle(z)}}$
provides the addition formula $\displaystyle{p=\frac{y \Delta(z)+y \Delta(z)}{1- \lambda_4 y^2 z^2}}$.
Using these functions, we obtain 
\begin{eqnarray}
\frac{1}{2}\left( \frac{\dd p}{\triangle(p)}+\frac{\dd q}{\triangle(q)}\right) =\frac{\dd y}{\triangle(y)}, \quad
\frac{1}{2}\left( \frac{\dd p}{\triangle(p)}-\frac{\dd q}{\triangle(q)}\right) =\frac{\dd z}{\triangle(z)} .
\label{6e3}
\end{eqnarray}
%
%%%%%%%%%%%%%%%%%%%%%%%%%%%%%%%%
\subsection{ $2 G(p,q) \pm \triangle(p) \triangle(q)$ as the function of  $y$ and $z$ }
Next, we calculate the necessary symmetric function of $p$ and $q$
\begin{eqnarray}
&&pq=\frac{y^2 \Delta(z)^2-z^2 \Delta(y)^2}{(1-y^2 z^2)^2}
=\frac{y^2(1-2 E z^2+z^4)-z^2(1-2 E y^2+y^4)}{(1-y^2 z^2)^2}
\nonumber\\
&&=\frac{y^2-z^2}{1-y^2 z^2}   ,
\label{6e4}\\
&&p^2+q^2=\frac{2\Big(y^2 \Delta(z)^2+z^2 \Delta(y)^2\Big)}{(1-y^2 z^2)^2}
=\frac{2\Big(y^2(1-2 E z^2+z^4)+z^2(1-2 E y^2+y^4)\Big)}{(1-y^2 z^2)^2}
\nonumber\\
&&=\frac{2\Big((1+y^2z^2)(y^2+z^2)-4E y^2 z^2\Big)}{1-y^2 z^2}  ,
\label{6e5}\\
&&\Delta(p) \Delta(q)=\mp \frac{(1+y^2z^2)^2-2E(1+y^2z^2) (y^2+z^2)+ (y^2+z^2)^2}{1-y^2 z^2}   .
\label{6e6}
\end{eqnarray}
For the ambiguity of sign in Eq.(\ref{6e6}), we take $(-1)$ sign.\footnote{G\"{o}pel takes $(+1)$ 
sign. In such case, in connection with the sign of Eq.(\ref{4e16}), we will see that 
the differential equation does not become of separable type.} \ Substituting the 
above expressions into
\begin{eqnarray}
2 G(p,q) \pm 2 \Delta(p) \Delta(q)
=2\Big(F(1+p^2 q^2)-C(p^2+q^2)+ 2Dpq\Big) \pm 2 \Delta(p) \Delta(q)  , 
\nonumber
\end{eqnarray}
we obtain $2 G(p,q) \pm 2 \Delta(p) \Delta(q)$ as the function of $y$ and $z$
\begin{eqnarray}
&&2 G(p,q) \pm 2 \Delta(p) \Delta(q)
\nonumber\\
&&= 
\frac{2(F \mp 1)}{(1-y^2 z^2)^2}\left( 1-2\frac{(C \mp E-D)y^2}{(F\mp 1)} +y^4 \right)
\left( 1-2\frac{(C \mp E+D)z^2}{(F\mp 1)} +z^4 \right)   ,
\label{6e7}
\end{eqnarray}
where we use the identity
$C^2-D^2+E^2+F^2-2 CEF=1$.
It is quite surprising that the $y$ dependence and the $z$ dependence 
becomes separable in the second and the third term of Eq.(\ref{6e7}).
Thus, we obtain $2 G(p,q) \pm 2 \Delta(p) \Delta(q)$ as the function of $y$ and $z$
\begin{eqnarray}
&&\sqrt{2 G(p,q) -2 \Delta(p) \Delta(q)}
\nonumber\\
&& 
=\frac{ \sqrt{2(F+1)} }{(1-y^2 z^2)} \sqrt{\left(1-2\frac{(C +E-D)y^2}{(F+1)} +y^4)
( 1-2\frac{(C + E+D)z^2}{(F+1)} +z^4\right)}
\nonumber\\
&&=\frac{\sqrt{2(F+1)}}{(1-y^2 z^2)}\sqrt{(1-2E_2y^2+y^4 ) ( 1-2E_3z^2+z^4)}   , 
\label{6e8}\\
&&\sqrt{2 G(p,q) +2 \Delta(p) \Delta(q)}
\nonumber\\
&& 
=\frac{ \sqrt{2(F-1)} }{(1-y^2 z^2)} \sqrt{\left(1-2\frac{(C -E-D)y^2}{(F-1)} +y^4)
( 1-2\frac{(C - E+D)z^2}{(F-1)} +z^4\right)}
\nonumber\\
&&=\frac{\sqrt{2(F-1)}}{(1-y^2 z^2)}\sqrt{(1-2E_1 y^2+y^4 ) ( 1-2E_4 z^2+z^4)}, 
\label{6e9}\\
&&
E_1=\frac{(C-E-D)}{(F-1)}, E_2=\frac{(C+E-D)}{(F+1)}, E_3=\frac{(C+E+D)}{(F+1)},
%\nonumber\\
%&&
E_4=\frac{(C-E+D)}{(F-1)},   
\nonumber
\end{eqnarray}
where we used the identity $C^2-D^2+E^2+F^2-2 C E F=1$ .

%%%%%%%%%%%%%%%%%%%%%%%%%
\subsection{The expression of $\varphi$ and $\psi$ with $y$ and $z$}
Next, we calculate $\varphi^2$ and $\psi^2$ as the function of $y$ and $z$. First, we obtain
\begin{eqnarray}
&&\varphi^2=(b+b_1)(1+p^2 q^2)-a (p^2+q^2)+2(c-c_1)pq
\nonumber\\
&&=(b+b_1) \left( 1+p^2q^2 -\frac{(C^2-1)}{(CF-E+D)} (p^2+ q^2)+\frac{2(CD-CE+F)}{(CF-E+D)}pq \right)
\nonumber\\
&&=\frac{b+b_1}{(1-y^2 z^2)^2}\left( 1-\frac{2(C-E-D)}{(F-1)} y^2+y^4 \right)
\left( 1-\frac{2(C+E+D)}{(F+1)} z^2+z^4 \right)
\nonumber\\
&&=\frac{b+b_1}{(1-y^2 z^2)^2}\left( 1-2 E_1 y^2+y^4 \right) 
\left( 1-2 E_3 z^2+z^4 \right)  .
\label{6e10}
\end{eqnarray}
Similarly,  we obtain
\begin{eqnarray}
&&\psi^2=(b-b_1)(1+p^2 q^2)-a (p^2+q^2)+2(c+c_1)pq
\nonumber\\
&&=(b-b_1) \left( 1+p^2q^2 -\frac{(C^2-1)}{(CF-E-D)} (p^2+ q^2)+\frac{2(CD+CE-F)}{(CF-E-D)}pq \right)
\nonumber\\
&&=\frac{b+b_1}{(1-y^2 z^2)^2}\left( 1-\frac{2(C+E-D)}{(F+1)} y^2+y^4 \right)
\left( 1-\frac{2(C-E+D)}{(F-1)} z^2+z^4 \right)
\nonumber\\
&&=\frac{b+b_1}{(1-y^2 z^2)^2}\left( 1-2 E_2 y^2+y^4 \right)
\left( 1-2 E_4 z^2+z^4 \right)   ,
\label{6e11}
\end{eqnarray}
where we used the identity $C^2-D^2+E^2+F^2-2 C E F=1$.
Thus, we obtain $\varphi$ and $\psi$ as the function of $y$ and $z$
\begin{eqnarray}
&& \varphi=\frac{\sqrt{b+b_1}}{(1-y^2 z^2)} 
\sqrt{( 1-2 E_1 y^2+y^4)( 1-2 E_3 z^2+z^4)}  ,
\label{6e12}\\
&& \psi=\frac{\sqrt{b-b_1}}{(1-y^2 z^2)} 
\sqrt{( 1-2 E_2 y^2+y^4)( 1-2 E_4 z^2+z^4)}  .
\label{6e13}
\end{eqnarray}
%

%%%%%%%%%%%%%%%%%%%%%%%%%
\subsection{The differential equation with $y$ and $z$ (Step III)}
Using Eqs.(\ref{6e3}), (\ref{6e8}), (\ref{6e9}), (\ref{6e12}) and (\ref{6e13}), 
we can express the necessary quantities with $y$ and $z$
\begin{eqnarray}
&& \frac{1}{2} \left(\frac{\dd p}{\triangle(p)}+\frac{\dd q}{\triangle(q)}\right)
=\frac{\dd y}{\triangle(y)}, \quad
\frac{1}{2} \left(\frac{\dd p}{\triangle(p)}-\frac{\dd q}{\triangle(q)}\right)
=\frac{\dd z}{\triangle(z)}   ,
\label{6e14}\\
&&\sqrt{2 G(p,q) -2 \Delta(p) \Delta(q)}
\nonumber\\
&&=
\frac{ \sqrt{2(F+1)} }{(1-y^2 z^2)} \sqrt{\left(1-2E_2 y^2+y^4\right)
\left( 1-2E_3 z^2+z^4\right)}  ,
\label{6e15}\\
&&\sqrt{2 G(p,q) +2 \Delta(p) \Delta(q)}
\nonumber\\
&&=\frac{ \sqrt{2(F-1)} }{(1-y^2 z^2)} \sqrt{\left(1-2 E_1y^2 +y^4\right)
\left( 1-2E_4 z^2+z^4\right)}  ,
\label{6e16}\\
&& \varphi=\frac{\sqrt{b+b_1}}{(1-y^2 z^2)} 
\sqrt{( 1-2 E_1 y^2+y^4)( 1-2 E_3 z^2+z^4)}   ,
\label{6e17}\\
&& \psi=\frac{\sqrt{b-b_1}}{(1-y^2 z^2)} 
\sqrt{( 1-2 E_2 y^2+y^4)( 1-2 E_4 z^2+z^4)}   .
\label{6e18}
\end{eqnarray}
Substituting these expressions into the following differential equations, 
\begin{eqnarray}
&& \frac{1}{\varphi} \left\{ \sqrt{2 G(p,q)-2 \triangle(p) \triangle(q)} 
\frac{1}{2}\left( \frac{\dd p}{\triangle(p)}+\frac{\dd q}{\triangle(q)}\right) \right.
\nonumber\\
&& \left. +\sqrt{2 G(p,q)+2 \triangle(p) \triangle(q)} 
\frac{1}{2}\left( \frac{\dd p}{\triangle(p)}-\frac{\dd q}{\triangle(q)}\right)\right\}
=\dd \mu, 
\label{6e19}\\
&& \frac{1}{\psi} \left\{ \sqrt{2 G(p,q)+2 \triangle(p) \triangle(q)} 
\frac{1}{2}\left( \frac{\dd p}{\triangle(p)}+\frac{\dd q}{\triangle(q)}\right) \right.
\nonumber\\
&& \left. +\sqrt{2 G(p,q)-2 \triangle(p) \triangle(q)} 
\frac{1}{2}\left( \frac{\dd p}{\triangle(p)}-\frac{\dd q}{\triangle(q)}\right)\right\}
=\dd \nu   , 
\label{6e20}
\end{eqnarray}
we obtain the separable differential equation of $y$ and $z$ in the form
\begin{eqnarray}
&&\dd \mu=\frac{1-y^2 z^2}{\sqrt{b+b_1}} \frac{1}{\sqrt{( 1-2 E_1 y^2+y^4)( 1-2 E_3 z^2+z^4)}}
\nonumber\\
&&\times \left( \frac{\sqrt{2(F+1)}}{1-y^2 z^2} \sqrt{(1-2E_2 y^2+y^4)
( 1-2E_3 z^2+z^4)} \frac{\dd y}{\triangle(y)}
\right.
\nonumber\\
&&\left. +\frac{\sqrt{2(F-1)}}{1-y^2 z^2} \sqrt{(1-2E_1 y^2+y^4)
( 1-2E_4 z^2+z^4)} \frac{\dd z}{\triangle(z)} \right)
\nonumber\\
&&=\frac{1}{\sqrt{b+b_1}} \left( \sqrt{\frac{2(F+1)(1-2 E_2 y^2+y^4)}{(1-2 E y^2+y^4)(1-2 E_1 y^2+y^4)}}
\dd y \right)
\nonumber\\
&&\left. +\sqrt{\frac{2(F-1)(1-2 E_4 z^2+z^4)}{(1-2 E z^2+z^4)(1-2 E_3 z^2+z^4)}} \dd z \right)  .
\label{6e21}
\end{eqnarray}
Similarly, we obtain 
\begin{eqnarray}
&&\dd \nu=\frac{(1-y^2 z^2)}{\sqrt{b-b_1}} \frac{1}{\sqrt{( 1-2 E_2 y^2+y^4)( 1-2 E_4 z^2+z^4)}}
\nonumber\\
&&\times \left( \frac{\sqrt{2(F-1)}}{(1-y^2 z^2)} \sqrt{(1-2E_1 y^2+y^4)
( 1-2E_4 z^2+z^4)} \frac{\dd y}{\triangle(y)}
\right.
\nonumber\\
&&\left. +\frac{\sqrt{2(F+1)}}{(1-y^2 z^2)} \sqrt{(1-2E_2 y^2+y^4)
( 1-2E_3 z^2+z^4)} \frac{\dd z}{\triangle(z)} \right)
\nonumber\\
&&=\frac{1}{\sqrt{b-b_1}} \left( \sqrt{\frac{2(F-1)(1-2 E_1 y^2+y^4)}{(1-2 E y^2+y^4)(1-2 E_2 y^2+y^4)}}
\dd y \right.
\nonumber\\
&&\left. +\sqrt{\frac{2(F+1)(1-2 E_3 z^2+z^4)}{(1-2 E z^2+z^4)(1-2 E_4 z^2+z^4)}} \dd z \right).
\label{6e22}
\end{eqnarray}
It is quite surprising that, if we use the functions $y$ and $z$, 
the differential equations become of the separable type.

%%%%%%%%%%%%%%%%%%%%%%%%%%%%%%%%%%%%%%%%%%%%%%%%%%%%%%%%%%%%%%%%%%%%%%%%%%%%%%
%%%%%%%%%%%%%%%%%%%%%%%%%%   Section 7  %%%%%%%%%%%%%%%%%%%%%%%%%%%%%%%%%%%%%%%%%
%%%%%%%%%%%%%%%%%%%%%%%%%%%%%%%%%%%%%%%%%%%%%%%%%%%%%%%%%%%%%%%%%%%%%%%%%%%%%
%%%%%%%%%%%%%%%%%%%%%%%%%%%%%%%%%%%%%%%%%%%%%%%%%%%%%%%%%%%%
\section{The differential equation with $y$ and $y'$ (Step IV)}
\setcounter{equation}{0}
Next, we change the function in order that it provide the same type of Abelian 
differential. 

%%%%%%%%%%%%%%%%%%%%%%%%%%%%
\subsection{The change of the function from $z$ to $y'$} 
The differential equation in the previous section provides
\begin{eqnarray}
&&\dd \mu=
\frac{1}{\sqrt{b+b_1} } \left( \frac{\sqrt{2(F+1)} (1-2 E_2 y^2+y^4)}
{\sqrt{(1-2 E y^2+y^4)(1-2 E_1 y^2+y^4)(1-2 E_2 y^2+y^4)}}
\dd y  \right.
\nonumber\\
&&
\left. +\frac{\sqrt{2(F-1)}(1-2 E_4 z^2+z^4)}
{\sqrt{(1-2 E z^2+z^4)(1-2 E_3 z^2+z^4)(1-2 E_4 z^2+z^4)}} \dd z \right)   ,
\label{7e1}\\
&&\dd \nu=
\frac{1}{\sqrt{b-b_1} } \left( \frac{\sqrt{2(F-1)} (1-2 E_1 y^2+y^4)}
{\sqrt{(1-2 E y^2+y^4)(1-2 E_1 y^2+y^4)(1-2 E_2 y^2+y^4)}}
\dd y  \right.
\nonumber\\
&&
\left. +\frac{\sqrt{2(F+1)}(1-2 E_3 z^2+z^4)}
{\sqrt{(1-2 E z^2+z^4)(1-2 E_3 z^2+z^4)(1-2 E_4 z^2+z^4)}} \dd z \right) .
\label{7e2}
\end{eqnarray}
For the function $y$, we obtain the Abelian differential of the type \\
$$\frac{\dd y}{\sqrt{(1-2E y^2+y^4)(1-2E_1 y^2+y^4)(1-2E_2 y^2+y^4)}}.$$
While, for the function $z$, we obtain the Abelian differential of the type 
$$\frac{\dd z}{\sqrt{(1-2E z^2+z^4)(1-2E_3 z^2+z^4)(1-2E_4 z^2+z^4)}}. $$
However, the type of the Abelian differential is different.

We keep $y$ to be the same function, yet we change the function
$z$ into $y'$, which cause the  M\"{o}bius transformation
of $z^2$, and make the Abel function of the same type.
We parametrize $2E$ with $\alpha$ in the form $2E=\alpha^2+1/\alpha^2$. 
For the constant $e$, we obtain $\Delta(e)=1-2 E e^2+e^4=(1-e^2/\alpha^2)
(1-e^2/\beta^2)$, $\alpha \beta=1$.
We change the function from $y$, $z$ to $y$, $y'$ ($y$ is not changed) by the following differential equation
\begin{eqnarray}
\frac{\dd z}{\triangle(z)}=\frac{\dd y'}{\triangle(y')}
+\frac{\dd e}{\triangle(e)}=\frac{\dd y'}{\triangle(y')}  , 
\label{7e3}
\end{eqnarray}
where we use $\dd e=0$ because $e$ is the constant.
This provides
\begin{eqnarray}
z=\frac{y' \triangle(e)+e \triangle(y')}{1-e^2 y'^2}  .
\label{7e4}
\end{eqnarray}
Next, we choose $e=\alpha$, thus we obtain $\triangle(e)=\triangle(\alpha)=0$, hence
$z$ is given as the function of $y'$
\begin{eqnarray}
z=\frac{\alpha \triangle(y')}{1-\alpha^2 y'^2}  .
\label{7e5}
\end{eqnarray}
Considering the square of $z$, we obtain
\begin{eqnarray}
z^2=\frac{\alpha^2-y'^2}{1-\alpha^2 y'^2}  .
\label{7e6}
\end{eqnarray}
Thus, $y'^2$ is the M\"{o}bius transformation of $z^2$.
Hence, we obtain
\begin{eqnarray}
\triangle(z)^2=(1-\frac{z^2}{\alpha^2}) (1-\frac{z^2}{\alpha^2}) 
=\left( \frac{(1-\alpha^2)y'}{\alpha (1-\alpha^2 y'^2)} \right)^2   ,
\label{7e7}
\end{eqnarray}
which provides $\Delta(z)$ as the function of $y$ in the following form 
\begin{eqnarray}
\triangle(z)=\frac{(1-\alpha^2)y'}{\alpha (1-\alpha^2 y'^2)} .
\label{7e8}
\end{eqnarray}
Using Eq.(\ref{7e5}) and Eq.(\ref{7e8}) , we obtain the necessary quantities as
the function of $y$ and $y'$ 
\begin{eqnarray}
&&p=\frac{ y \triangle(z)+z \triangle(y)}{1-y^2 z^2}
=\frac{\beta(\alpha^4-1)y y'+\alpha \triangle(y)  \triangle(y')}
{1-\alpha^2(y^2+y'^2)+y^2 y'^2},
\label{7e9}\\
&&q=\frac{ y \triangle(z)-z \triangle(y)}{1-y^2 z^2}
=\frac{\beta(\alpha^4-1)y y'-\alpha \triangle(y)  \triangle(y')}
{1-\alpha^2(y^2+y'^2)+y^2 y'^2},
\label{7e10}\\
&&pq=\frac{ y^2-z^2}{1-y^2 z^2}
=\frac{y^2+y'^2-\alpha^2(1+y^2 y'^2)}
{1-\alpha^2(y^2+y'^2)+y^2 y'^2}  .
\label{7e11}
\end{eqnarray}
We will use these relations in the next section.

%%%%%%%%%%%%%%%%%%%%%%%%%%%%%%%%%%%%%%%%%%%%
\subsection{The differential equation with $y$ and $y'$ (Step IV)}
We calculate the following quantity as the function of $y'$
\begin{eqnarray}
\frac{1}{2}(z^2+\frac{1}{z^2})=
\frac{1}{2}\left (\frac{\alpha^2-y'^2}{1-\alpha^2 y'^2}+
\frac{1-\alpha^2 y'^2}{\alpha^2-y'^2}\right)
=\frac{E-2 y'^2 +E y'^4}{1-2 E y'^2 +y'^4}   .
\nonumber
\end{eqnarray}
Thus, we obtain 
\begin{eqnarray}
\frac{1-2 E_3 z^2+z^4}{2z^2}=\frac{1}{2}(z^2+\frac{1}{z^2})-E_3
=\frac{(E-E_3)-2(1-E E_3) y'^2 +(E-E_3) y'^4}{1-2 E y'^2 +y'^4}  .
\label{7e12}
\end{eqnarray}
Similarly, we obtain
\begin{eqnarray}
\frac{1-2 E_4 z^2+z^4}{2z^2}=\frac{1}{2}(z^2+\frac{1}{z^2})-E_4
=\frac{(E-E_4)-2(1-E E_4) y'^2 +(E-E_4) y'^4}{1-2 E y'^2 +y'^4}  .
\label{7e13}
\end{eqnarray}
By using the identity
\begin{eqnarray}
1-E(E_1+E_3)+E_1 E_3=0, \quad
1-E(E_2+E_4)+E_2 E_4=0, \quad
\frac{E-E_4}{E-E_3}=\frac{F+1}{F-1}  ,
\nonumber
\end{eqnarray}
we obtain in the form
\begin{eqnarray}
&&\frac{1-2 E_3 z^2+z^4}{2z^2}=\frac{(E-E_3)-2(1-E E_3) y'^2 +(E-E_3) y'^4}{1-2 E y'^2 +y'^4}
\nonumber\\
&&=\frac{(E-E_3)(1-2 E_1 y'^2 +y'^4)}{1-2 E y'^2 +y'^4},
\label{7e14}\\
&&\frac{1-2 E_4 z^2+z^4}{2z^2}=\frac{(E-E_4)-2(1-E E_4) y'^2 +(E-E_4) y'^4}{1-2 E y'^2 +y'^4}
\nonumber\\ 
&&=\frac{(E-E_4)(1-2 E_2 y'^2 + y'^4)}{1-2 E y'^2 +y'^4} .
\label{7e15}
\end{eqnarray}
Taking the ratio of Eq.(\ref{7e14}) and Eq.(\ref{7e15}), we obtain the expression
\begin{eqnarray}
&&\frac{1-2 E_3 z^2+z^4}{1-2 E_4 z^2+z^4}
=\left(\frac{F-1}{F+1}\right)
\left(\frac{1-2 E_1 y'^2+y'^4}{1-2 E_2 y'^2+y'^4}\right) ,
\label{7e16}\\
&&\frac{1-2 E_4 z^2+z^4}{1-2 E_3 z^2+z^4}
=\left(\frac{F+1}{F-1}\right)
\left(\frac{1-2 E_2 y'^2+y'^4}{1-2 E_1 y'^2+y'^4}\right)  .
\label{7e17}
\end{eqnarray}
Therefore, we obtain the differential equation of the same type Abelian differential  
\begin{eqnarray}
&&\dd \mu=\frac{1}{\sqrt{b+b_1}} 
\left( \sqrt{\frac{2(F+1)(1-2 E_2 y^2+y^4)}{(1-2 E y^2+y^4)(1-2 E_1 y^2+y^4)}}\dd y
+\sqrt{\frac{2(F-1)(1-2 E_4 z^2+z^4)}{(1-2 E z^2+z^4)(1-2 E_3 z^2+z^4)}} \dd z \right)
\nonumber\\
&&=\sqrt{\frac{2(F+1)}{b+b_1}} 
\left( \sqrt{\frac{(1-2 E_2 y^2+y^4)}{(1-2 E y^2+y^4)(1-2 E_1 y^2+y^4)}}\dd y \right.
\nonumber\\
&&\left. +\sqrt{\frac{(1-2 E_2 y'^2+y'^4)}{(1-2 E y'^2+y'^4)(1-2 E_1 y'^2+y'^4)}} \dd y' \right)  .
\label{7e18}
\end{eqnarray}
Similarly, we obtain 
\begin{eqnarray}
&&\dd \nu=\frac{1}{\sqrt{b-b_1}} 
\left( \sqrt{\frac{2(F-1)(1-2 E_1 y^2+y^4)}{(1-2 E y^2+y^4)(1-2 E_2 y^2+y^4)}} \dd y
+\sqrt{\frac{2(F+1)(1-2 E_3 z^2+z^4)}{(1-2 E z^2+z^4)(1-2 E_4 z^2+z^4)}} \dd z \right)
\nonumber\\
&&=\sqrt{\frac{2(F-1)}{b-b_1}}
 \left( \sqrt{\frac{(1-2 E_1 y^2+y^4)}{(1-2 E y^2+y^4)(1-2 E_2 y^2+y^4)}} \dd y \right.
\nonumber\\
&&\left. +\sqrt{\frac{(1-2 E_1 y'^2+y'^4)}{(1-2 E y'^2+y'^4)(1-2 E_2 y'^2+y'^4)}} \dd y' \right)  .
\label{7e19}
\end{eqnarray}
Thus, for $\dd y$ and $\dd y'$, we have the same type Abelian differential
of the form
$$\frac{\dd x}{\sqrt{(1-2 E x^2+x^4)(1-2 E_1 x^2+x^4)(1-2 E_2 x^2+x^4)}}.$$

%%%%%%%%%%%%%%%%%%%%%%%%%%%%%%%%%%%%%%%%%%%%%%%%%%%%%%%%%%%%%%%%%%%%%%%%%%%%%%
%%%%%%%%%%%%%%%%%%%%%%%%%%   Section 8  %%%%%%%%%%%%%%%%%%%%%%%%%%%%%%%%%%%%%%%%%
%%%%%%%%%%%%%%%%%%%%%%%%%%%%%%%%%%%%%%%%%%%%%%%%%%%%%%%%%%%%%%%%%%%%%%%%%%%%%
%%%%%%%%%%%%%%%%%%%%%%%%%%%%%%%%%%%%%%%%%%%%%%%%%%%%%%%%%%%%
\section{The differential equation with $x$ and $x'$(Step V) }
\setcounter{equation}{0}
%%%%%%%%%%%%%%%%%%%%%%%%%%%
\subsection{The differential equation of the genus two Jacobi's inversion problem}
By changing the functions $y$, $y'$ into $x$, $x'$ in such a way 
as $\dd x$ and $\dd x'$ becomes the same Abelian differential of the fifth-degree hyperelliptic curve.
For that purpose, we change the function in the form
\begin{eqnarray}
x=\left( \frac{1-y^2}{1+y^2}\right)^2, \quad x'=\left( \frac{1-y'^2}{1+y'^2}\right)^2   . 
 \label{8e1}
\end{eqnarray}
We denote in the following 
\begin{eqnarray}
m=\frac{E+1}{E-1}, \quad, m_1=\frac{E_1+1}{E_1-1}, \quad, m_2=\frac{E_2+1}{E_2-1}   .
\label{8e2}
\end{eqnarray}
Thus, we obtain the simple expression
\begin{eqnarray}
\frac{1-m_2 x}{1-m_1 x}=\left(\frac{E_1-1}{E_2-1}\right) 
\left(\frac{1-2 E_2 y^2+y^4}{1-2 E_1 y^2+y^4}\right) , 
\frac{1-m_1 x}{1-m_2 x}=\left(\frac{E_2-1}{E_1-1}\right)  
\left(\frac{1-2 E_1 y^2+y^4}{1-2 E_2 y^2+y^4}\right)  .
\label{8e3}
\end{eqnarray}
Therefore, we obtain the simple expression 
\begin{eqnarray}
\sqrt{ \frac{1-2 E_1 y^2+y^4}{1-2 E_2 y^2+y^4} }=
\sqrt{ \frac{E_1-1}{E_2-1} } \sqrt{ \frac{1-m_1 x}{1-m_2 x} }, \ 
\sqrt{ \frac{1-2 E_1 y'^2+y'^4}{1-2 E_2 y'^2+y'^4} }=
\sqrt{ \frac{E_1-1}{E_2-1} } \sqrt{ \frac{1-m_1 x'}{1-m_2 x'} },
\label{8e4}
\end{eqnarray}
and
\begin{eqnarray}
\sqrt{ \frac{1-2 E_2 y^2+y^4}{1-2 E_1 y^2+y^4} }=
\sqrt{ \frac{E_2-1}{E_1-1} } \sqrt{ \frac{1-m_2 x}{1-m_1 x} }, \ 
\sqrt{ \frac{1-2 E_2 y'^2+y'^4}{1-2 E_1 y'^2+y'^4} }=
\sqrt{ \frac{E_2-1}{E_1-1} } \sqrt{ \frac{1-m_2 x'}{1-m_1 x'} }   .
\label{8e5}
\end{eqnarray}
By differentiating  $y$, we obtain
\begin{eqnarray}
&&\dd x=-\frac{8((1-y^2)y \dd y}{(1+y^2)^3}, \quad \sqrt{x}=\frac{1-y^2}{1+y^2}, \quad
\sqrt{1-x}=\frac{2y}{1+y^2}, 
\nonumber\\
&&\sqrt{1-mx}
=\sqrt{ \frac{-2}{E-1} } \frac{\sqrt{1-2E y^2+y^4}}{1+y^2}    . 
\label{8e6}
\end{eqnarray}
Combining these relations, we obtain the connection of the differential of $\dd x$ and $\dd y$
\begin{eqnarray}
\frac{\dd x}{\sqrt{x(1-x)(1-mx)}}=\frac{2\sqrt{2}\sqrt{-1}
\sqrt{E-1}\ \dd y}{\sqrt{1-2 E y^2+y^4}}   .
 \label{8e7}
\end{eqnarray}
Similarly, we obtain 
\begin{eqnarray}
\frac{\dd x'}{\sqrt{x'(1-x')(1-mx')}}=\frac{2\sqrt{2}\sqrt{-1}
\sqrt{E-1}\ \dd y'}{\sqrt{1-2 E y'^2+y'^4}}  .
 \label{8e8}
\end{eqnarray}
Using Eqs.(\ref{7e18}), (\ref{7e19}), (\ref{8e7}) and (\ref{8e8}),  we finally obtain
the simplified differential equation
\begin{eqnarray}
&&\hskip -20mm \sqrt{\frac{2(F+1)}{b+b_1}} 
\left( \sqrt{\frac{(1-2 E_2 y^2+y^4)}{(1-2 E y^2+y^4)(1-2 E_1 y^2+y^4)}}\dd y
+\sqrt{\frac{(1-2 E_2 y'^2+y'^4)}{(1-2 E y'^2+y'^4)(1-2 E_1 y'^2+y'^4)}} \dd y' \right)
\nonumber\\
&&\hskip -20mm =\frac{1}{2\sqrt{b+b_1}}  \sqrt{ \frac{(F+1)(1-E_2)}{(1-E)(1-E_1)} }
\left(  \sqrt{  \frac{(1-m_2 x)}{x(1-x)(1-mx)}  } \dd x +\sqrt{  \frac{(1-m_2 x')}{x'(1-x')(1-mx')}  } \dd x'  \right)
=\dd \mu, 
\label{8e9}\\
&&\hskip -20mm \sqrt{\frac{2(F-1)}{b-b_1}}
 \left( \sqrt{\frac{(1-2 E_1 y^2+y^4)}{(1-2 E y^2+y^4)(1-2 E_2 y^2+y^4)}}\dd y
+\sqrt{\frac{(1-2 E_1 y'^2+y'^4)}{(1-2 E y'^2+y'^4)(1-2 E_2 y'^2+y'^4)}} \dd y' \right)
\nonumber\\
&&\hskip -20mm =\frac{1}{2\sqrt{b+b_1}}  \sqrt{ \frac{(F+1)(1-E_2)}{(1-E)(1-E_1)} }
\left(  \sqrt{  \frac{(1-m_2 x)}{x(1-x)(1-mx)}  } \dd x +\sqrt{  \frac{(1-m_2 x')}{x'(1-x')(1-mx')}  } \dd x'  \right)
=\dd \nu  .
\label{8e10}
\end{eqnarray}
Therefore, we obtain the differential equation of the genus two Jacobi's inversion problem  
\begin{eqnarray}
&&
\frac{(1-m_2 x) \dd x}{\sqrt{f_5(x)}}
+\frac{(1-m_2 x' ) \dd x'}{\sqrt{f_5(x')}}
=2 \sqrt{b+b_1} \sqrt{\frac{(1-E)(1-E_1)}{(F+1)(1-E_2)}} \dd \mu=\dd \widehat{\mu}, 
\label{8e11}\\
&&
\frac{(1-m_1 x) \dd x}{\sqrt{f_5(x)}}
+\frac{(1-m_1 x' ) \dd x'}{\sqrt{f_5(x')}}
=2 \sqrt{b-b_1} \sqrt{\frac{(1-E)(1-E_2)}{(F-1)(1-E_1)}} \dd \nu=\dd \widehat{\nu}    , 
\label{8e12}\\
&&
f_(x)=x(1-x)(1-mx)(1-m_1 x) (1-m_2 x)    .
\label{8e13}
\end{eqnarray} 
%

%%%%%%%%%%%%%%%%%%%%%%%%%%%%%%
\subsection{The expression of $x$ and $x'$ as the function of $p$ and $q$}
Next, we express $x$, $x'$ with $p$, $q$.
From Eqs.(\ref{7e9})-(\ref{7e11}), we obtain the necessary quantities as the function of $y$ and $y'$
\begin{eqnarray}
p+q= \frac{2 \beta (\alpha^4-1) y y'}{ 1-\alpha^2 (y^2+y'^2)+y^2 y'^2}, \quad
pq= \frac{y^2+y'^2-\alpha^2(1+y^2 y'^2)}{ 1-\alpha^2 (y^2+y'^2)+y^2 y'^2},
\nonumber
\end{eqnarray}
which provides
\begin{eqnarray}
1+\alpha^2 pq=\frac{(1-\alpha^4)(1+y^2 y'^2)}{1+y^2 y'^2-\alpha^2(y^2+y'^2)}, \quad
1+\beta^2 pq=\frac{(\beta^2-\alpha^2)(y^2+ y'^2)}{1+y^2 y'^2-\alpha^2(y^2+y'^2)}  . 
\label{8e14}
\end{eqnarray}
Thus, we obtain
\begin{eqnarray}
\frac{p+q}{1+\alpha^2 pq}=-\frac{2\beta y y'}{(1+y^2 y'^2)}, \quad
\frac{p+q}{1+\beta^2 pq}=-\frac{2\alpha y y'}{(y^2+ y'^2)}   .
\label{8e15}
\end{eqnarray}
Hence, we can, in principle, express $y^2+y'^2$ and $y y'$ as the function of $p$ and $q$,
which implies that $x=(1-y^2)^2/(1+y^2)^2$ and $x'=(1-y'^2)^2/(1+y'^2)^2$ can 
be expressed as the function of $p$ and $q$. 
For our purpose, we make the following combination
\begin{eqnarray}
f=\left(\frac{1-y y'}{1+y y'}\right)^2, \quad
g=\left(\frac{y- y'}{y+ y'}\right)^2,
\label{8e16} 
\end{eqnarray}
which provides the following nice factorization property
\begin{eqnarray}
f=\frac{(1+\alpha p)(1+\alpha q)}{(1-\alpha p)(1-\alpha q)}, \quad
g=\frac{(1+\beta p)(1+\beta q)}{(1-\beta p)(1-\beta q)}  .
\label{8e17}
\end{eqnarray}
Later, we use $f+g$, $1+fg$, $\sqrt{fg}$. Thus, we calculate these quantities
\begin{eqnarray}
&&f+g=\frac{2( p^2q^2+(\alpha-\beta)^2 p q -(p^2+q^2)+1)}
{(1-\alpha p)(1-\alpha q)(1-\beta p)(1-\beta q)}
\nonumber\\
&&=\frac{2( p^2 q^2+(2E-2) p q -(p^2+q^2)+1)}
{(1-\alpha p)(1-\alpha q)(1-\beta p)(1-\beta q)},  
\label{8e18}\\
&&1+fg=\frac{2( p^2q^2+(\alpha+\beta)^2 p q +(p^2+q^2)+1)}
{(1-\alpha p)(1-\alpha q)(1-\beta p)(1-\beta q)}
\nonumber\\
&&=\frac{2( p^2q^2+(2 E+2) p q +(p^2+q^2)+1)}
{(1-\alpha p)(1-\alpha q)(1-\beta p)(1-\beta q)}, 
\label{8e19}\\
&&\sqrt{fg}=\pm \frac{\Delta(p) \Delta(q)}
{(1-\alpha p)(1-\alpha q)(1-\beta p)(1-\beta q)}  ,
\label{8e20}
\end{eqnarray}
where we used $\alpha^2+\beta^2=2 E$, $\alpha \beta=1$.
We take $(+1)$ sign in Eq.(\ref{8e20}).\footnote{G\"{o}pel take $(-1)$ sign in Eq.(\ref{8e20}). } 
\ Using Eq.(\ref{8e16}), we obtain  
\begin{eqnarray}
\left(\frac{1-y y'}{1+y y'}\right)=\sqrt{f}, \quad
\left(\frac{y- y'}{y+ y'}\right)=\sqrt{g}, 
\label{8e21}
\end{eqnarray}
which provides $y y'$, $y/y'$ as the function of $\sqrt{f}$, $\sqrt{g}$
\begin{eqnarray}
yy'=\frac{1-\sqrt{f}}{1+\sqrt{f}}, \quad
\frac{y}{y'}=\frac{1+\sqrt{g}}{1-\sqrt{g}}.
\label{8e22}
\end{eqnarray}
Multiplying and dividing the first and the second term of Eq.(\ref{8e22}), we obtain
\begin{eqnarray}
y^2=\frac{(1-\sqrt{f})(1+\sqrt{g})}{(1+\sqrt{f})(1-\sqrt{g})} , \quad
y'^2=\frac{(1-\sqrt{f})(1-\sqrt{g})}{(1+\sqrt{f})(1+\sqrt{g})} .
\label{8e23}
\end{eqnarray}
Thus, we obtain the desired form 
\begin{eqnarray}
\frac{1-y^2}{1+y^2}=\frac{\sqrt{f}-\sqrt{g}}{1-\sqrt{fg}} , \quad
\frac{1-y'^2}{1+y'^2}=\frac{\sqrt{f}+\sqrt{g}}{1+\sqrt{fg}} .
\label{8e24}
\end{eqnarray}
Therefore, $x$ and $x'$ are expressed as the function of $p$ and $q$
\begin{eqnarray}
&&x=\left(\frac{1-y^2}{1+y^2}\right)^2=\frac{f+g-2\sqrt{fg}}{1+fg-2\sqrt{fg}}
=\frac{E(1+p^2 q^2)-(p^2+q^2)-\Delta(p) \Delta(q)}{(E+1)(pq+1)^2}  ,
\label{8e25}\\
&&
x'=\left(\frac{1-y'^2}{1+y'^2}\right)^2=\frac{f+g+2\sqrt{fg}}{1+fg+2\sqrt{fg}}
=\frac{E(1+p^2 q^2)-(p^2+q^2)+\Delta(p) \Delta(q)}{(E+1)(pq+1)^2}  .
\label{8e26}
\end{eqnarray}
%

%%%%%%%%%%%%%%%%%%%%%%%%%%%
\subsection{The solution of the genus two Jacobi's inversion problem}
From the differential equation
\begin{eqnarray}
\frac{(1-m_2 x) \dd x}{\sqrt{f_5(x)}} +\frac{(1-m_2 x' ) \dd x'}{\sqrt{f_5(x')}}
=\dd \widehat{\mu}, 
\quad 
\frac{(1-m_1 x) \dd x}{\sqrt{f_5(x)}}+\frac{(1-m_1 x' ) \dd x'}{\sqrt{f_5(x')}}
=\dd \widehat{\nu}  , 
\nonumber
\end{eqnarray} 
we can rearrange in the standard Jacobi's inversion problem of the form
\begin{eqnarray}
\frac{\dd x}{\sqrt{f_5(x)}} +\frac{\dd x'}{\sqrt{f_5(x')}}
=\dd U_1, \quad 
\frac{x \dd x}{\sqrt{f_5(x)}} + \frac{x' \dd x'}{\sqrt{f_5(x')}}
=\dd U_2   , 
\label{8e27}
\end{eqnarray} 
where we obtain the expression $U_1=k_1 u_1+ k_2 u_2$, $U_2=k_3 u_1+ k_4 u_2$. We inversely express
$u_1=\ell_1 U_1+ \ell_2 U_2$, $u_2=\ell_3 U_1+ \ell_4 U_2$. 
Thus, the solution of the Jacobi's inversion problem is given by
\begin{eqnarray}
&&\wp_{22}(U_1,U_2)=x+x'
=\frac{2(E(1+p^2q^2)-(p^2+q^2))}{(E+1)(pq+1)^2}
\nonumber\\
&&=\frac{2E}{(E+1)} -\frac{2 (p^2+q^2-2 E pq)}{(E+1)(1+pq)^2}, 
\label{8e28}\\
&&\wp_{12}(U_1,U_2)=-x x'
=-\frac{(E-1)}{(E+1)} +\frac{4(E-1) pq}{(E+1)(1+pq)^2}, 
\label{8e29}
\end{eqnarray}
where
\begin{eqnarray}
&&p=\frac{S_1(u_1,u_2)}{P_1(u_1,u_2)}
=\frac{\vartheta[\begin{array}{cc} 1 \ 1 \\ 0\ \ 1 \\ \end{array}](\ell_1 U_1+ \ell_2 U_2,\ell_3 U_1+ \ell_4 U_2)}
{\vartheta[\begin{array}{cc} 0 \ 0 \\ 0\ \ 1 \\ \end{array}](\ell_1 U_1+ \ell_2 U_2, \ell_3 U_1+ \ell_4 U_2)}, 
\label{8e30}\\
&&q=\frac{S_2(u_1,u_2)}{P_2(u_1,u_2)}
=\frac{\vartheta[\begin{array}{cc} 1 \ 1 \\ 1\ \ 0 \\ \end{array}](\ell_1 U_1+ \ell_2 U_2, \ell_3 U_1+ \ell_4 U_2)}
{\vartheta[\begin{array}{cc} 0 \ 0 \\ 1\ \ 0 \\ \end{array}](\ell_1 U_1+ \ell_2 U_2,\ell_3 U_1+ \ell_4 U_2)} .
\label{8e31}
\end{eqnarray}
For the explicit functional form for $\wp_{22}(U_1, U_2), \wp_{12}(U_1, U_2)$,  we must carefully determine
$\ell_i, (i=1,2,3,4)$ by the $u_1=0, u_2=0$ value of various theta functions and derivative 
of the various theta functions.

%%%%%%%%%%%%%%%%%%%%%%%%%%%%%%%%%%%%%%%%%%%%%%%%%%%%%%%%%%%%
%%%%%%%%%%%%%%%%%%%%%%%%%%%%%%%%%%%%%%%%%%%%%%%%%%%%%%%%%%%%%%%%%%%%%%%%%%%%%%
%%%%%%%%%%%%%%%%%%%%%%%%%%   Section 9  %%%%%%%%%%%%%%%%%%%%%%%%%%%%%%%%%%%%%%%%%
%%%%%%%%%%%%%%%%%%%%%%%%%%%%%%%%%%%%%%%%%%%%%%%%%%%%%%%%%%%%%%%%%%%%%%%%%%%%%
\section{Summary and Discussions}
\setcounter{equation}{0}
%%%%%%%%%%%%%%%%%%%%%%%%%%%%%%%%%%%%%%%%%%%%%%%%%%%%%%%%%
 In the previous paper, we reviewed the Rosenhain's paper
to the Jacobi's inversion problem for the genus two 
hyperelliptic integral.
In this paper, we have reviewed the G\"{o}pel's paper 
to the Jacobi's inversion problem for the genus two 
hyperelliptic integral.

In the Rosenhain's approach,  
the  Riemann's  addition formula of the hyperelliptic theta function is used.  
The key identity of the Rosenhain's approach is the 
three quadratic theta identities.
In the G\"{o}pel's approach, 
the addition formula by the duplication method is used.
The key identity of the G\"{o}pel's approach  is the three quartic theta identities,  
Eq.(\ref{3e25}), Eq.(\ref{5e11}), Eq.(\ref{5e20}).  One of these identities, Eq.(\ref{3e25}),  
is the quartic Kummer surface relation. Three quartic identities are given 
in the form 
\begin{eqnarray}
&&P_1^4+S_1^4+P_2^4+S_2^4 -2F(P_1^2 P_2^2+S_1^2 S_2^2)+2C(P_1^2 S_2^2+P_2^2 S_1^2)
\nonumber\\
&&-2E(P_1^2 S_1^2+P_2^2 S_2^2)-4D P_1 S_1 P_2 S_2=0,
\nonumber\\
&&P_3^2 S_3^2 +P_0^2 S_0^2-b(P_1^2 P_2^2+S_1^2 S_2^2)+a(P_1^2 S_2^2+P_2^2 S_1^2)
-2c P_1 S_1 P_2 S_2=0,
\nonumber\\
&&2 P_3 S_3 P_0 S_0-b_1 (P_1^2 P_2^2+S_1^2 S_2^2)+2c_1P_1 S_1 P_2 S_2=0.
\nonumber
\end{eqnarray}

Starting from the genus two hyperelliptic theta functions, G\"{o}pel takes 
several steps to obtain the differential equation of the genus two Jacobi's 
inversion equation.
In step I, by using the addition formula of the genus two hyperelliptic theta
functions, the derivative formula is obtained; thus the starting differential 
equation is obtained. In step II, by using three quartic theta identities, 
the differential equation of only $p$ and $q$ is obtained.
In step III, by using the addition formula of the elliptic function 
of the form 
$\displaystyle{\frac{\dd p}{\Delta(p)} =\frac{\dd y}{\Delta(y)}+\frac{\dd z}{\Delta(z)}}$
and 
$\displaystyle{\frac{\dd q}{\Delta(q)} =\frac{\dd y}{\Delta(y)}-\frac{\dd z}{\Delta(z)}}$
, the separable type differential equation of $y$ and $z$ is obtained.
In step IV, by using the addition formula of the elliptic function 
of the form 
$\displaystyle{\frac{\dd z}{\Delta(z)} 
=\frac{\dd y'}{\Delta(y')}+\frac{\dd e}{\Delta(e)}=\frac{\dd y'}{\Delta(y')}}$ with 
some special constant $e$, the differential equation of $y$ and $y'$ with 
the same type Abelian differential is obtained.
In step V, via $\displaystyle{x=\left( \frac{1-y^2}{1+y^2} \right)^2 }$,
$\displaystyle{x'=\left( \frac{1-y'^2}{1+y'^2} \right)^2 }$, the differential 
equation of $x$ and $x'$ with the Abelian differential of  the fifth-degree hyperelliptic 
curve is obtained, i.e., the differential equation of the genus 
two Jacobi's inversion problem is obtained.

In the paper on the comments of the G\"{o}pel's paper, Jacobi pointed out 
the issue of the number of modules~\cite{Jacobi5}.
For the genus $g$ hyperelliptic curve 
$y^2=\lambda_{2g+1} x^{2g+1}+\lambda_{2g} x^{2g}+\cdots+\lambda_1 x+\lambda_0$,
we have the $2g+2$ constant coefficients. However, we have the freedom to 
change, i) constant scale of $x$, ii) constant shift of $x$, iii) take the ratio of the coefficients, 
in order to provide the standard form.
For example, in the $g=1$ case, i) provides 
$y^2=4 x^3+\lambda_2 x^2+\lambda_1 x+\lambda_0$.
While, ii) provide $y^2=4 x^3-g_2 x-g_3$. In the step iii),  $g_2^3/(g_2^3-27 g_3^2)$
provides $\tau=\omega_3/\omega_1$ from $\omega_3$ and $\omega_1$
by taking the ratio of the constant coefficients.
Thus, the number of the modules from the hyperelliptic curve is 
$2g+2-3=2g-1$.  While, in the $g$-variable hyperelliptic theta function,
in general,  there are modules of the form $\tau_{i j}(=\tau_{j i}), (i, j=1, 2,\cdots , g)$; thus, the number 
of the modules is $g(g+1)/2$.
Hence,  for $g \ge 3$, 
the number of the module for the hyperelliptic function $g(g+1)/2$  
becomes greater than that of the hyperelliptic curve $2g-1$.

\vskip 20mm
%%%%%%%%%%%%%%%%%%%%%%%%%%%%%%%%%%%%%%%%%%%%%%%%%%%%
%\newpage
%%%%%%%%%%%%%%%%%%%%%%%%% References %%%%%%%%%%%%%%%%%%%

%
%%%%%%%%%%%%%%%%%%%%%%%%%%%%%%%%%%%%%%%%%%%%%%%
%%%%%%%%%%%%%%%%%%%%%%%%%%%%%%%%%%%%%%%%%%%%%%%
%%%%%%%%%%%%%%%%%%%%%%%%%%%%%%%%%%%%%%%%%%%%%%%
%%%%%%%%%%%%%%%%%%%%%%%%%%%%%%%%%%%%%%%%%%%%%%%
%%%%%%%%%%%%%%%%%%%%%%%%%%%%%%%%%%%%%%%%%%%%%%%
%%%%%%%%%%%%%%%%%%%%%%%%%%%%%%%%%%%%%%%%%%%%%%%

\end{document}